\documentclass[10pt,reqno]{amsart}

\usepackage[
  letterpaper,
  left=1in,
  right=1in,
  top=1in,
  bottom=1in
]{geometry}

\usepackage{amsmath,amssymb,amsfonts,amsthm}
\usepackage{mathtools}
\usepackage{bm}

\usepackage{graphicx}
\graphicspath{{figures/}}
\usepackage{booktabs}
\usepackage{enumitem}
\usepackage{caption}
\usepackage{subcaption}
\usepackage{algorithm}
\usepackage{algpseudocode}
\usepackage{float}
\usepackage{placeins}
\usepackage{tikz}
\usetikzlibrary{arrows.meta,positioning,calc,decorations.pathreplacing}
\usepackage{pgfplots}
\pgfplotsset{compat=1.18}
\usepgfplotslibrary{fillbetween}
\definecolor{accent}{HTML}{8C1C1C}  

\usepackage[hidelinks]{hyperref}

\newtheorem{theorem}{Theorem}[section]
\newtheorem{proposition}[theorem]{Proposition}
\newtheorem{lemma}[theorem]{Lemma}
\newtheorem{corollary}[theorem]{Corollary}
\theoremstyle{definition}

\newtheorem{assumption}[theorem]{Assumption}
\theoremstyle{remark}
\newtheorem{remark}[theorem]{Remark}

\newcommand{\RR}{\mathbb{R}}
\newcommand{\PP}{\mathbb{P}}
\newcommand{\EE}{\mathbb{E}}

\newcommand{\cS}{\mathcal{S}}
\newcommand{\cE}{\mathcal{E}}

\newcommand{\Phibar}{\bar\Phi}
\newcommand{\G}[1]{\ifcsname gen@#1\endcsname\csname gen@#1\endcsname\else
  \textbf{??#1??}\PackageWarning{generated}{Missing generated value #1}\fi}
\expandafter\def\csname gen@gauss-Rsim\endcsname{40}
\expandafter\def\csname gen@gauss-Rtune\endcsname{10}
\expandafter\def\csname gen@gauss-Tsim\endcsname{250000}
\expandafter\def\csname gen@gauss-age-after\endcsname{0.094}
\expandafter\def\csname gen@gauss-age-before\endcsname{0.87}
\expandafter\def\csname gen@gauss-age-m\endcsname{32}
\expandafter\def\csname gen@gauss-age-p\endcsname{0.999}
\expandafter\def\csname gen@gauss-asym-ratio-hi\endcsname{0.72}
\expandafter\def\csname gen@gauss-asym-ratio-lo\endcsname{0.66}
\expandafter\def\csname gen@gauss-collapse-d2\endcsname{0.0158}
\expandafter\def\csname gen@gauss-collapse-hi\endcsname{0.0176}
\expandafter\def\csname gen@gauss-collapse-k\endcsname{10^{-3}}
\expandafter\def\csname gen@gauss-collapse-lo\endcsname{0.0172}
\expandafter\def\csname gen@gauss-collapse-m-hi\endcsname{29}
\expandafter\def\csname gen@gauss-collapse-m-lo\endcsname{28}
\expandafter\def\csname gen@gauss-filt-ratio-k4\endcsname{0.94}
\expandafter\def\csname gen@gauss-g2-excess-hi\endcsname{18}
\expandafter\def\csname gen@gauss-g2-excess-lo\endcsname{10}
\expandafter\def\csname gen@gauss-grid\endcsname{3000}
\expandafter\def\csname gen@gauss-grid-dev\endcsname{0.006}
\expandafter\def\csname gen@gauss-ill-g2arg\endcsname{22}
\expandafter\def\csname gen@gauss-ill-mA\endcsname{14}
\expandafter\def\csname gen@gauss-ill-p\endcsname{0.99}
\expandafter\def\csname gen@gauss-kappa-hi\endcsname{10^{-1}}
\expandafter\def\csname gen@gauss-kappa-lo\endcsname{10^{-4}}
\expandafter\def\csname gen@gauss-nconfig\endcsname{28}
\expandafter\def\csname gen@gauss-nm-A\endcsname{0.096}
\expandafter\def\csname gen@gauss-nm-R\endcsname{40}
\expandafter\def\csname gen@gauss-nm-T\endcsname{250000}
\expandafter\def\csname gen@gauss-nm-best-bur\endcsname{19}
\expandafter\def\csname gen@gauss-nm-bur\endcsname{0.073}
\expandafter\def\csname gen@gauss-nm-geo\endcsname{0.087}
\expandafter\def\csname gen@gauss-nm-mA\endcsname{14}
\expandafter\def\csname gen@gauss-nm-mean\endcsname{100}
\expandafter\def\csname gen@gauss-nm-ratio-hi\endcsname{1.28}
\expandafter\def\csname gen@gauss-nm-ratio-lo\endcsname{1.18}
\expandafter\def\csname gen@gauss-nm-ratio-se-max\endcsname{0.002}
\expandafter\def\csname gen@gauss-nm-reg\endcsname{0.092}
\expandafter\def\csname gen@gauss-nm-se-max\endcsname{0.0007}
\expandafter\def\csname gen@gauss-offline-hi\endcsname{3.0}
\expandafter\def\csname gen@gauss-offline-lo\endcsname{1.6}
\expandafter\def\csname gen@gauss-robust-double\endcsname{24}
\expandafter\def\csname gen@gauss-robust-half\endcsname{10}
\expandafter\def\csname gen@gauss-rule-excess-max\endcsname{1.4}
\expandafter\def\csname gen@gauss-sandwich-gap-max\endcsname{2.9}
\expandafter\def\csname gen@gauss-sim-ratio-e-hi\endcsname{1.11}
\expandafter\def\csname gen@gauss-sim-ratio-e-lo\endcsname{1.03}
\expandafter\def\csname gen@gauss-sim-ratio-e-se-max\endcsname{0.01}
\expandafter\def\csname gen@gauss-sim-ratio-u-hi\endcsname{1.38}
\expandafter\def\csname gen@gauss-sim-ratio-u-lo\endcsname{1.15}
\expandafter\def\csname gen@gauss-sim-se-max\endcsname{0.0003}
\expandafter\def\csname gen@gauss-skew-mstar-dev\endcsname{12}
\expandafter\def\csname gen@gauss-skew-pstar-dev\endcsname{3.2}
\expandafter\def\csname gen@gauss-skew-ratio\endcsname{3}
\expandafter\def\csname gen@gauss-smoother-gain-hi\endcsname{2.0}
\expandafter\def\csname gen@gauss-smoother-gain-lo\endcsname{1.3}
\expandafter\def\csname gen@kreg-Rsim\endcsname{40}
\expandafter\def\csname gen@kreg-Tsim\endcsname{250000}
\expandafter\def\csname gen@kreg-asym-hi\endcsname{0.73}
\expandafter\def\csname gen@kreg-asym-lo\endcsname{0.72}
\expandafter\def\csname gen@kreg-bd-pstar\endcsname{0.095}
\expandafter\def\csname gen@kreg-cost-adj-max\endcsname{0.507}
\expandafter\def\csname gen@kreg-cost-jump-hi\endcsname{0.771}
\expandafter\def\csname gen@kreg-cost-jump-lo\endcsname{0.755}
\expandafter\def\csname gen@kreg-em-40\endcsname{1.1\cdot10^{-3}}
\expandafter\def\csname gen@kreg-filt-bd-k3\endcsname{0.0133}
\expandafter\def\csname gen@kreg-filt-uni-k3\endcsname{0.0111}
\expandafter\def\csname gen@kreg-jump-excess\endcsname{13}
\expandafter\def\csname gen@kreg-kappa\endcsname{0.01}
\expandafter\def\csname gen@kreg-kappa-hi\endcsname{0.03}
\expandafter\def\csname gen@kreg-kappa-lo\endcsname{10^{-4}}
\expandafter\def\csname gen@kreg-kz-hi\endcsname{1.07}
\expandafter\def\csname gen@kreg-kz-lo\endcsname{0.94}
\expandafter\def\csname gen@kreg-kz-se-max\endcsname{0.01}
\expandafter\def\csname gen@kreg-nconfig\endcsname{18}
\expandafter\def\csname gen@kreg-rev-fwd\endcsname{0.0890}
\expandafter\def\csname gen@kreg-rev-limit\endcsname{1.10}
\expandafter\def\csname gen@kreg-rev-ratio\endcsname{1.063}
\expandafter\def\csname gen@kreg-rev-rev\endcsname{0.0837}
\expandafter\def\csname gen@kreg-rule-A-excess\endcsname{0.5}
\expandafter\def\csname gen@kreg-two-law-excess-hi\endcsname{0.4}
\expandafter\def\csname gen@kreg-uni-pstar\endcsname{0.107}
\expandafter\def\csname gen@kreg-win-bd-k3\endcsname{0.0182}
\expandafter\def\csname gen@kreg-win-uni-k3\endcsname{0.0207}
\expandafter\def\csname gen@kreg-winratio-bd-k3\endcsname{1.37}
\expandafter\def\csname gen@kreg-winratio-cyc-k3\endcsname{2.14}
\expandafter\def\csname gen@kreg-winratio-se-max\endcsname{0.02}
\expandafter\def\csname gen@kreg-winratio-uni-k3\endcsname{1.87}
\expandafter\def\csname gen@pm1-A-excess-max\endcsname{13.6}
\expandafter\def\csname gen@pm1-A-excess-med\endcsname{0.0}
\expandafter\def\csname gen@pm1-R\endcsname{120}
\expandafter\def\csname gen@pm1-T\endcsname{30000}
\expandafter\def\csname gen@pm1-e-filt-hi\endcsname{1.041}
\expandafter\def\csname gen@pm1-e-filt-lo\endcsname{1.010}
\expandafter\def\csname gen@pm1-g2-excess-hi\endcsname{73}
\expandafter\def\csname gen@pm1-maxdev\endcsname{0.002}
\expandafter\def\csname gen@pm1-maxz\endcsname{2.8}
\expandafter\def\csname gen@pm1-mstar-hi\endcsname{17}
\expandafter\def\csname gen@pm1-mstar-lo\endcsname{1}
\expandafter\def\csname gen@pm1-nA-exact\endcsname{18}
\expandafter\def\csname gen@pm1-nagree-sim\endcsname{30}
\expandafter\def\csname gen@pm1-nconfig\endcsname{30}
\expandafter\def\csname gen@pm1-nones\endcsname{4}
\expandafter\def\csname gen@pm1-tie-delta\endcsname{0.80}
\expandafter\def\csname gen@pm1-tie-prob\endcsname{0.42}
\expandafter\def\csname gen@pm1-u-filt-hi\endcsname{1.20}
\expandafter\def\csname gen@pm1-u-filt-lo\endcsname{1.10}
\expandafter\def\csname gen@refl-M-hi\endcsname{160}
\expandafter\def\csname gen@refl-M-lo\endcsname{20}
\expandafter\def\csname gen@refl-R\endcsname{60}
\expandafter\def\csname gen@refl-T\endcsname{30000}
\expandafter\def\csname gen@refl-b-best-hi\endcsname{0.111}
\expandafter\def\csname gen@refl-b-best-lo\endcsname{0.109}
\expandafter\def\csname gen@refl-b-pred-hi\endcsname{0.133}
\expandafter\def\csname gen@refl-b-raw-hi\endcsname{0.276}
\expandafter\def\csname gen@refl-b-wall-hi\endcsname{21}
\expandafter\def\csname gen@refl-b-wall-lo\endcsname{4}
\expandafter\def\csname gen@refl-delta\endcsname{1.2}
\expandafter\def\csname gen@refl-err5-hi\endcsname{0.288}
\expandafter\def\csname gen@refl-err5-lo\endcsname{0.147}
\expandafter\def\csname gen@refl-err5-proxy\endcsname{0.108}
\expandafter\def\csname gen@refl-m5\endcsname{5}
\expandafter\def\csname gen@refl-minbound-lo\endcsname{0.55}
\expandafter\def\csname gen@refl-mstar-list\endcsname{13, 7, 7 and 7}
\expandafter\def\csname gen@refl-p\endcsname{0.97}
\expandafter\def\csname gen@refl-pa-maxdev\endcsname{0.002}
\expandafter\def\csname gen@refl-ph5-hi\endcsname{0.58}
\expandafter\def\csname gen@refl-ph5-lo\endcsname{0.13}
\expandafter\def\csname gen@refl-proxy-best\endcsname{0.108}
\expandafter\def\csname gen@refl-proxy-mstar\endcsname{5}
\expandafter\def\csname gen@refl-se-max\endcsname{0.0017}
\expandafter\def\csname gen@refl-wall-hi\endcsname{41}
\expandafter\def\csname gen@refl-wall-lo\endcsname{9}
{generated/numbers.tex not found; run code/run_all.py}}

\newcommand{\authorcontact}[3]{%
  \par\addvspace{\bigskipamount}
  \noindent\begin{minipage}{\linewidth}
  #1\par
  \textsc{#2}\par
  \if\relax\detokenize{#3}\relax\else\textit{Email address:}\space\texttt{#3}\fi
  \end{minipage}\par
}

\title[Sliding-Window Detection of Hidden Regimes]{Sliding-Window Detection of Hidden Regimes under Stochastic Switching}

\author{A. Pambukyan}

\author{H. Nersisyan}

\date{\today}

\begin{document}

\begin{abstract}
A sliding-window average is a simple estimator of the current state of a hidden regime-switching process, such as a hidden Markov model. A longer window averages out more noise but is more likely to contain a switch, and the question is what such a switch costs. We analyse the decoder that assigns the average of the last \(m\) observations to the nearest regime mean. For two regimes with conditionally independent, stochastically ordered observations, a switch inside the window, averaged uniformly over its positions and directions, causes an error with probability exactly \(1/2\). Consequently, for every stationary switching law, the error is at most \(\bar P(m)+\bar s(m-1)/2\), where \(\bar P(m)\) is the stationary-weighted single-regime error and \(\bar s\) the switch rate. Under Markov switching, time reversal shows that a window whose oldest and current regimes differ is misclassified with probability exactly \(1/2\), however many switches it contains, and the bound is attained up to a factor \(p_{\min}^{m-1}\). For Gaussian noise with common variance, minimizing the bound gives an explicit window length which, in units of \(\sigma^2/\Delta^2\), depends only on the expected number of switches in that time. As switching probabilities vanish proportionally, the optimal error is \(\bar s\log(1/\bar s)/(2I)\,(1+o(1))\), with \(I\) the single-regime Chernoff exponent. In the Gaussian case this is twice the Bayes error. For symmetric switching, exponentially weighted moving averages reduce this factor to at most \(2\log2\approx1.39\), and two-block weights adapted to the switch rate attain the Bayes rate. For three or more regimes, long-window switch costs are set by the geometry of the regime means. For Markov-modulated random walks with \(\pm1\) steps on a reflecting lattice, a position-aware decoder learns as fast as without walls under holding boundaries, while steps at bouncing boundaries carry none.

\end{abstract}

\maketitle

\section{Introduction}\label{sec:intro}

Many processes alternate between a few regimes that are not observed directly. An unobserved regime \((E_t)\) with finitely many values determines the law of each observation \(Y_t\), and when \((E_t)\) is a Markov chain the pair is a hidden Markov model. Asset returns alternate between bull and bear markets \cite{Hamilton1989,AngTimmermann2012}, animals between foraging and travelling \cite{LangrockEtAl2012}, tracked molecules between diffusive states \cite{PerssonEtAl2013}, and queues between periods of high and low load \cite{PrabhuZhu1989}. Often the current regime has to be estimated from noisy observations. The Bayes filter does this optimally, but it requires the transition matrix and the full law of the observations \cite{CappeMoulinesRyden2005}. In practice simpler estimators, such as a sliding window or a discounted average of recent observations, can be useful, as in sliding-window and discounted bandit algorithms \cite{GarivierMoulines2011,GajaneOrtnerAuer2018} and in rolling-window regime detection \cite{SheikhBhat2026}. They, however, require choosing how much of the past to use. Once its length is fixed, a window decoder uses only the regime means. A longer window averages out more noise but is more likely to contain a switch, and bounding the two effects separately (Proposition~\ref{thm:recoverability}) treats every window that contains a switch as an error. Such a bound determines neither the optimal window nor its error. What is less obvious is how much a switch inside the window actually costs.

A window whose \(k\) most recent observations come from the current regime and whose other \(m-k\) come from the previous one has the same composition as the mirrored window, in which the two regimes exchange roles, and the decoder errs on complementary events in the two. Averaged uniformly over its \(m-1\) possible positions and both directions, a switch inside the window therefore causes an error with probability exactly one half, whatever the emission laws. This is a probabilistic counterpart, in the presence of noise, of the lag of \((m-1)/2\) steps of a moving average. Under Markov switching, time reversal gives an exact conditional statement. Any decoder unchanged by reversing the window misclassifies a window whose oldest and current regimes differ with probability one half, however many switches it contains. The asymptotics of the Bayes filter are due to Khasminskii and Zeitouni \cite{KhasminskiiZeitouni1996}, and the baseline bound is standard. To our knowledge, the exact cost of a switch and its consequences, summarized below, have not been stated before.

The question arose for Markov-modulated random walks (MMRWs), in which a hidden finite-state Markov chain \((E_t)\), the environment, selects at each step the transition kernel of an observed walk \((X_t)\). Without reflection, such a walk is a Markov additive process \cite{Asmussen2003} and a hidden Markov model observed through its increments \cite{CappeMoulinesRyden2005}. Pambukyan, Saudagar and Kumar \cite{PambukyanSaudagarKumar2026} established ergodicity of finite MMRWs on a reflecting lattice and observed that standard supervised learners recover the regime from short windows of increments. Sections~\ref{sec:half-window} and~\ref{sec:k-regimes} apply to any regime-switching process with conditionally independent observations. The reflected MMRW, in which the law of an increment also depends on the position of the walk, is treated as an application in Section~\ref{sec:reflected}.

\subsection{Main results}

A useful dimensionless parameter is the expected number of switches per resolution time
\(\kappa=\bar s\sigma^2/\Delta^2\), where \(\bar s\) is the probability of a switch at a given step and \(\sigma^2/\Delta^2\) is the number of observations needed to separate the regime means.

\emph{Two regimes} (Section~\ref{sec:half-window}). Suppose that, given the regime path, the observations are independent and their two laws are stochastically ordered. Then the error probability of the uniform window of length \(m\) is at most
\[
A(m)=\bar P(m)+\bar s\,\frac{m-1}{2},
\]
where \(\bar P(m)=\sum_i\pi_iP^{(i)}(m)\) is the stationary-weighted single-regime error, \(P^{(i)}(m)\) being the error of a window that lies entirely in regime~\(i\). The bound holds for every stationary switching law. Under Markov switching the error is also at least \(p_{\min}^{m-1}A(m)\), where \(p_{\min}\) is the smallest diagonal entry of the transition matrix, and both bounds hold without any ordering of the emission laws. The factor \(p_{\min}^{m-1}\) tends to one at the optimal window as switches become rare (Theorem~\ref{thm:half-window}). Minimizing \(A\) gives a near-optimal window, explicit for Gaussian noise with common variance. As all switching probabilities tend to zero at a fixed ratio, and under exponential-moment conditions on the emission laws, the optimal error is
\[
P^\star=\frac{\bar s\log(1/\bar s)}{2I}\,(1+o(1)),
\]
where \(I\) is the Chernoff exponent of a window without switches (Corollary~\ref{cor:window-rule}). For Gaussian noise with common variance \(P^\star=4\kappa\log(1/\kappa)(1+o(1))\), twice the error of the Bayes filter \cite{KhasminskiiZeitouni1996}. The factor two is due to uniform weighting. Exponential weights reduce it to at most \(2\log2\approx1.39\) (Proposition~\ref{prop:shape}), weights whose shape changes with \(\kappa\) attain the Bayes rate in the limit (Proposition~\ref{prop:two-block}), and for nearly noiseless \(\pm1\) steps the factor tends to one. For an agent who acts on the decoded regime, these factors are also ratios of avoidable regret (Remark~\ref{rem:regret}).

\emph{More than two regimes} (Section~\ref{sec:k-regimes}). A window that contains a switch may now be assigned to a third regime. Defining the cost of a switch as the error averaged over its position in the window, the error is bounded by the single-regime error plus the expected total cost of the switches in the window, up to the probability of two or more switches (Theorem~\ref{thm:switch-cost}). Asymptotically, a switch between adjacent regimes costs \(1/2\), and a switch from \(\mu_j\) to \(\mu_i\) costs \(1-h_{ij}/|\mu_i-\mu_j|\), where \(h_{ij}\) is the half-gap of \(\mu_i\) on the side of \(\mu_j\) (Proposition~\ref{prop:switch-geometry}). The window must be long enough to separate the closest regimes, whereas the Bayes filter detects a large jump faster than a small one, so the ratio of the two errors grows with the size of the jumps. For Gaussian noise with common variance and equally spaced means, a switch across \(d\) gaps costs the window \(2d(2d-1)\) times what it costs the filter.

\emph{Application to the reflected walk} (Section~\ref{sec:reflected}). Near a reflecting boundary the law of an increment depends on the position of the walk, and visits to the boundary can reduce the accuracy of a window of raw increments. A decoder that weights each increment by its starting position is governed by an information clock, which advances at each step by the signal-to-noise ratio of that step (Proposition~\ref{prop:position-aware}). For \(\pm1\) steps and holding boundaries, every step has the same signal-to-noise ratio as in the interior, because a blocked step is observed as a zero increment. Under bouncing boundaries the steps at a boundary carry none.

Section~\ref{sec:experiments} checks the results against exact computation and simulation. Proofs and additional experiments are in the appendices.

\subsection{Related work}

\emph{Filtering of hidden Markov chains.} Khasminskii and Zeitouni \cite{KhasminskiiZeitouni1996} obtained the asymptotic error of the Bayes filter for rare switches and any irreducible chain (see also \cite{Chigansky2006}). Their upper bound uses simple filters that sum log-likelihood ratios over windows of length of order \(\log(1/\varepsilon)/\mathrm{KL}\) and need the emission laws but not the transition rates. We quantify what is lost by using only the regime means. Adaptive social learning, in which agents discount past log-likelihoods geometrically, attains the order \(\varepsilon\log(1/\varepsilon)\) of the Bayes error under rare transitions \cite{KhammassiEtAl2025}. For a single agent, two symmetric regimes and Gaussian noise with common variance it reduces to exponential weights, and Proposition~\ref{prop:shape} gives its constant for adaptation times of order at most the mean sojourn, which include those used there. This constant exceeds that of the Bayes filter by a factor \(2\log2\). For two hidden states, Abraham, Gassiat and Naulet \cite{AbrahamGassiatNaulet2023} prove matching minimax bounds for parameter estimation, and Gassiat, Kaddouri and Naulet \cite{GassiatKaddouriNaulet2025} relate the risk of clustering the hidden states to the Bayes risk of classifying them. Both concern the unsupervised problem, in which the emission laws are unknown. Here the regime means are given, and the object of study is the decoder.

\emph{Window methods.} Moving averages and weighted sums are classical tools for detecting a change in quality control \cite{Lai1974,BassevilleNikiforov1993}, with detection delay and false alarms as criteria. In the Window-Limited CUSUM of Xie, Moustakides and Xie \cite{XieMoustakidesXie2023}, a sliding window estimates the unknown post-change parameter inside a test for a single change, and its length is chosen for first-order optimality as the average run length tends to infinity. Moving sums also estimate multiple change points, including random ones generated by hidden Markov and regime-switching models \cite{EichingerKirch2018} and switches in the drift of multivariate processes \cite{KirchKlein2023}. The bandwidth balances detection power against resolution for change points with a minimum spacing, and the criteria are consistency and localization. Schwier, Brooks and Griffin \cite{SchwierBrooksGriffin2011} bound the window needed to distinguish switching Markov models, trading the false positives of short windows against the delay of long ones. In non-stationary bandits, sliding windows and discounting are tuned to the number of changes \cite{GarivierMoulines2011}, and Boot and Pick \cite{BootPick2018} weight observations optimally for forecasting under Markov switching. In this paper the regime means are known and switches recur, the window is the decoder itself, and the criterion is the stationary probability of misclassifying the current regime, windows with several switches included. Under Markov switching with Gaussian noise of common variance, the continuous window rule, in units of \(\sigma^2/\Delta^2\), is determined by \(\kappa\). Sheikh and Bhat \cite{SheikhBhat2026} detect a predictive regime from a rolling window, with a bound that counts every switch in the window and a window length logarithmic in the switching frequency. Under the assumptions of Theorem~\ref{thm:half-window}, the switching term is half as large, and the optimal window and its asymptotic constant are explicit.

\emph{Decision-making under hidden regimes.} Regime recovery is a component of decision-making in hidden-mode Markov decision processes \cite{ChoiYeungZhang2000}, regime-switching bandits \cite{ZhouXiongChenGao2021} and portfolio trading under regime switching \cite{DaiZhangZhu2010}.

\section{Model and Decoders}\label{sec:model}

The latent environment is a time-homogeneous Markov chain \((E_t)_{t\ge 0}\) on \(\cE:=\{1,\dots,K\}\) with transition matrix \(Q\) and stationary distribution \(\pi\). Unless stated otherwise, \(Q\) is irreducible and aperiodic. The regime \(E_t\) governs the law of a real observation \(Y_t\). In the motivating example, the MMRW, the observations are the increments
\[
Y_t:=X_t-X_{t-1}
\]
of a walk \((X_t)_{t\ge 0}\) on the lattice \(\cS:=\{0,1,\dots,M\}\) with reflecting boundaries, in which the regime \(E_t\) selects the kernel of the step from \(X_{t-1}\) to \(X_t\). The joint process \((X_t,E_t)\) is then a finite-state Markov chain on \(\cS\times\cE\), ergodic under mild conditions \cite{PambukyanSaudagarKumar2026}. Sections~\ref{sec:half-window} and~\ref{sec:k-regimes} use only the structure below and apply to any regime-switching process of this form, and Section~\ref{sec:reflected} returns to the walk.

\subsection{Regime-conditioned observations}

Each regime \(i\) has a mean \(\mu_i\in\RR\), and
\begin{equation}
Y_t=\mu_{E_t}+\xi_t .
\label{eq:emission-model}
\end{equation}

\begin{assumption}[Autonomous regimes and sub-Gaussian residuals]
\label{ass:residual-model}
\leavevmode
\begin{enumerate}[label=(\roman*)]
\item The environment is autonomous: for every \(t\ge 0\) and \(j\in\cE\),
\(\PP(E_{t+1}=j\mid E_0,\ldots,E_t,Y_1,\ldots,Y_t)=Q_{E_tj}\).
\item With \(\mathcal H_t:=\sigma(E_1,\ldots,E_t,Y_1,\ldots,Y_{t-1})\), for every \(t\ge 1\) and \(\lambda\in\RR\),
\[
\EE[\xi_t\mid\mathcal H_t]=0,
\qquad
\EE\!\left[e^{\lambda\xi_t}\mid\mathcal H_t\right]\le \exp\!\left(\frac{\sigma^2\lambda^2}{2}\right)\quad\text{a.s.}
\]
\end{enumerate}
\end{assumption}

Part~(ii) makes \((\xi_t)\) a conditionally \(\sigma^2\)-sub-Gaussian martingale-difference sequence, so the residuals may depend on past observations and on the regime path so far. It holds when the residuals are independent given the regimes with a \(\sigma^2\)-sub-Gaussian law depending only on the current regime, as in the observation models of Section~\ref{sec:experiments}. The separation between regimes is
\begin{equation}
\Delta:=\min_{i\neq j}|\mu_i-\mu_j|>0.
\label{eq:Delta}
\end{equation}

\subsection{Persistence and switches per resolution time}

Local persistence is measured by \(p_{\min}:=\min_i Q_{ii}\) and, in stationarity, by the switch rate
\begin{equation}
\bar s:=\PP_\pi(E_{t+1}\neq E_t)=\sum_{i\in\cE}\pi_i(1-Q_{ii}).
\label{eq:switch-rate}
\end{equation}
A stationary chain is regarded as defined for all \(t\in\mathbb Z\), which is possible by Kolmogorov's extension theorem. This matters only for events that reach before time \(0\), such as the age of the current regime exceeding \(t\).

An average of \(m\) residuals is sub-Gaussian with variance proxy \(\sigma^2/m\), so its standard deviation is at most \(\sigma/\sqrt m\), which equals the separation \(\Delta\) when \(m=\sigma^2/\Delta^2\). We call \(\sigma^2/\Delta^2\) the resolution time and
\begin{equation}
\kappa:=\bar s\,\frac{\sigma^2}{\Delta^2}
\label{eq:kappa}
\end{equation}
the expected number of switches per resolution time.

\subsection{Window decoders}

Given a window length \(m\ge 1\), the window average and the nearest-mean decoder are
\begin{equation}
\bar Y_{t,m}:=\frac1m\sum_{s=t-m+1}^{t}Y_s,
\qquad
\widehat E_t\in\arg\min_{i\in\cE}|\bar Y_{t,m}-\mu_i|,
\label{eq:decoder}
\end{equation}
with ties broken at random. For two regimes, labelled so that \(\mu_1<\mu_2\), we also use general linear statistics. Let \(\mathcal T\) be a finite set of times, \(w=(w_s)_{s\in\mathcal T}\) nonnegative weights with \(\sum_s w_s=1\), and \(\theta:=(\mu_1+\mu_2)/2\). The linear decoder with weights \(w\) is
\begin{equation}
S_t:=\sum_{s\in\mathcal T}w_sY_s,
\qquad
\widehat E_t:=
\begin{cases}
2,& S_t>\theta,\\
1,& S_t<\theta,
\end{cases}
\label{eq:linear-decoder}
\end{equation}
with ties broken by a fair coin. It contains the trailing uniform window \(\mathcal T=\{t-m+1,\ldots,t\}\), \(w_s=1/m\), which is \eqref{eq:decoder} for two regimes, exponential weights \(w_{t-j}\propto\beta^j\), and, offline, the centred window \(\mathcal T=\{t-h,\ldots,t+h\}\).

The event that the window contains no switch is
\begin{equation}
\mathrm{NoSwitch}_{t,m}
:=
\{E_{t-m+1}=\cdots=E_t\}.
\label{eq:noswitch-def}
\end{equation}

\begin{remark}[Known regime means]
\label{rem:known-means}
Sections~\ref{sec:half-window}--\ref{sec:reflected} treat \(\mu\) (and, for window design, \(\sigma\) and the switching probabilities) as known, which isolates the recovery mechanism. Section~\ref{sec:plug-in} shows that errors in the estimated switch rate cost the window rule little. How much of one trajectory is needed to estimate these quantities is a separate question, which we do not treat here.
\end{remark}

Table~\ref{tab:assumptions} lists the assumptions of each result.

\begin{table}[t]
\centering
\small
\caption{Assumptions of the main results. \emph{Conditionally independent}: given the regime path, the observations are independent with \(Y_s\sim F_{E_s}\). \emph{Rare switches}: the switching probabilities tend to zero with fixed ratios.}
\label{tab:assumptions}
\begin{tabular}{@{}llll@{}}
\toprule
Result & Switching & Observations & Validity\\
\midrule
Proposition~\ref{thm:recoverability} & Markov & sub-Gaussian residuals (Assumption~\ref{ass:residual-model}) & every \(m\)\\
Theorem~\ref{thm:half-window}(i) & stationary, two regimes & conditionally independent, ordered & every \(m\)\\
Theorem~\ref{thm:half-window}(ii) & Markov, two regimes & conditionally independent, unordered & every \(m\)\\
Corollary~\ref{cor:window-rule}(iii) & Markov, two regimes & as above, exponential moments & rare switches\\
Proposition~\ref{prop:shape} & symmetric Markov & Assumption~\ref{ass:residual-model}, Gaussian for equality & rare switches\\
Proposition~\ref{prop:two-block} & symmetric Markov & Assumption~\ref{ass:residual-model} & rare switches\\
Proposition~\ref{prop:centred} & symmetric Markov & Assumption~\ref{ass:residual-model} & every \(m\)\\
Theorem~\ref{thm:switch-cost} & stationary, Markov in (ii) & conditionally independent & every \(m\)\\
Proposition~\ref{prop:switch-geometry} & -- & finite variances, Gaussian for the expansion & \(m\to\infty\)\\
Corollary~\ref{cor:k-window}(ii) & Markov & conditionally independent, exponential moments & rare switches\\
Propositions~\ref{prop:reflected},~\ref{prop:position-aware} & Markov & increments of the reflected walk & every \(m\)\\
\bottomrule
\end{tabular}
\end{table}

\section{Two Regimes: The Half-Window Law}\label{sec:half-window}

For two regimes, the error of a window decoder obeys a simple law that, under Markov switching, is asymptotically sharp at the optimal window. We begin with the usual bound.

\subsection{A baseline bound}

By the Markov property, a window that starts in regime~\(i\) contains no switch with probability \(Q_{ii}^{\,m-1}\ge p_{\min}^{\,m-1}\), so that under stationarity
\begin{equation}
\PP(\mathrm{NoSwitch}_{t,m}^c)=1-\sum_{i\in\cE}\pi_iQ_{ii}^{\,m-1},
\label{eq:noswitch-exact}
\end{equation}
which equals \(1-p^{m-1}\) for the symmetric two-state chain \(Q(p)\).

\begin{proposition}[Baseline bound]
\label{thm:recoverability}
Under Assumption~\ref{ass:residual-model}, for every \(t\ge m\ge1\),
\begin{equation}
\PP(\widehat E_t\neq E_t)
\le
2\exp\!\left(-\frac{m\Delta^2}{8\sigma^2}\right)
+
\PP(\mathrm{NoSwitch}_{t,m}^c)
\le
2\exp\!\left(-\frac{m\Delta^2}{8\sigma^2}\right)
+1-p_{\min}^{\,m-1}.
\label{eq:recoverability-main}
\end{equation}
\end{proposition}

The bound holds for any number of regimes, and its proof, a Chernoff bound for the martingale of residuals, is in Appendix~\ref{app:proofs-local}. For two regimes only one tail of the residual sum leads to an error, and conditioning on the regime path as in Lemma~\ref{lem:path-mgf} shows that the factor~\(2\) can be dropped:
\begin{equation}
g_2(m):=\exp\!\left(-\frac{m\Delta^2}{8\sigma^2}\right)+1-p_{\min}^{\,m-1}.
\label{eq:g2}
\end{equation}
The first term bounds the error of a window without switches, and the second bounds the conditional error of a window that contains a switch by one. For two regimes, a switch averaged uniformly over its positions and both directions causes an error with probability one half.

\subsection{The half-window law}

From now on in this section \(K=2\), \(\mu_1<\mu_2\) and \(\theta:=(\mu_1+\mu_2)/2\). Write \(\tau_t:=\sup\{k\ge1:E_{t-k+1}=\cdots=E_t\}\in\{1,2,\ldots\}\cup\{\infty\}\) for the age of the current regime and \(\Phibar:=1-\Phi\) for the standard normal tail. When the observations are independent given the regimes, the half-window effect holds with equality when averaged uniformly over the positions and directions of a switch, whatever the emission laws.

\begin{theorem}[Half-window law]
\label{thm:half-window}
Let \((E_t)\) be a stationary process on two regimes, with \(\pi_i:=\PP(E_t=i)\) and switch rate \(\bar s:=\PP(E_{t+1}\neq E_t)\). Suppose that, given the whole regime path, the observations are independent with \(Y_s\sim F_{E_s}\), where \(F_1\) and \(F_2\) have means \(\mu_1<\mu_2\). For the trailing uniform window of length \(m\le t\), with ties broken by a fair coin, let \(P^{(i)}(m)\) be the error on a window that lies entirely in regime~\(i\), \(\bar P(m):=\pi_1P^{(1)}(m)+\pi_2P^{(2)}(m)\), and
\begin{equation}
A(m):=\bar P(m)+\bar s\,\frac{m-1}{2}.
\label{eq:A}
\end{equation}
\begin{enumerate}[label=(\roman*)]
\item If \(F_1\) is stochastically smaller than \(F_2\), then \(\PP(\widehat E_t\neq E_t)\le A(m)\).
\item If \((E_t)\) is a Markov chain with transition matrix \(Q\), and \(\lambda:=Q_{11}+Q_{22}-1\), then
\begin{equation}
\PP\bigl(\widehat E_t\neq E_t,\ E_{t-m+1}\neq E_t\bigr)=\tfrac12\,\PP\bigl(E_{t-m+1}\neq E_t\bigr)=\pi_1\pi_2\bigl(1-\lambda^{m-1}\bigr),
\label{eq:reversal}
\end{equation}
and consequently
\[
p_{\min}^{\,m-1}A(m)\le\sum_i\pi_iQ_{ii}^{\,m-1}P^{(i)}(m)+\pi_1\pi_2\bigl(1-\lambda^{m-1}\bigr)\le\PP(\widehat E_t\neq E_t)\le A(m),
\]
without any ordering of \(F_1\) and \(F_2\).
\end{enumerate}
\end{theorem}

\begin{proof}
For \(0\le k\le m\) and \(i\in\{1,2\}\), let \(\psi_i(k)\) be the error probability when regime~\(i\) is current, \(k\) of the \(m\) observations come from regime~\(i\) and \(m-k\) from the other regime, so that \(\psi_i(m)=P^{(i)}(m)\). Two properties of \(\psi_1\) and \(\psi_2\) drive the proof.

\emph{Monotonicity.} Under the stochastic-ordering assumption of part~(i), if regime~1 is current, the decoder errs with probability \(\EE[\mathbf 1\{S>m\theta\}+\tfrac12\mathbf 1\{S=m\theta\}]\), where \(S=\sum_sY_s\) is stochastically increasing in the number of observations from \(F_2\), and symmetrically for regime~2. Hence \(\psi_1\) and \(\psi_2\) are nonincreasing.

\emph{Pairing.} A window in which regime~1 is current with composition \((k,m-k)\) and a window in which regime~2 is current with composition \((m-k,k)\) contain the same laws, namely \(k\) draws from \(F_1\) and \(m-k\) from \(F_2\). The sum \(S\) has the same law in both, and the decoder errs on complementary events, so
\[
\psi_1(k)+\psi_2(m-k)=\PP(S>m\theta)+\PP(S<m\theta)+\PP(S=m\theta)=1,
\qquad\text{hence}\qquad
\sum_{k=1}^{m-1}\bigl(\psi_1(k)+\psi_2(k)\bigr)=m-1 .
\]

\emph{Upper bound.} On \(\{\tau_t\ge m,\ E_t=i\}\) the error is \(P^{(i)}(m)\), and this event has probability at most \(\pi_i\). On \(\{\tau_t=k,\ E_t=i\}\) with \(k<m\), at least \(k\) observations come from regime~\(i\), so by monotonicity the conditional error is at most \(\psi_i(k)\). This event requires a switch into regime~\(i\) at time \(t-k+1\), and by stationarity
\[
\PP(E_{s-1}=2,\ E_s=1)-\PP(E_{s-1}=1,\ E_s=2)=\PP(E_s=1)-\PP(E_{s-1}=1)=0,
\]
so a switch into either regime has probability \(\bar s/2\) at any time. Hence
\[
\PP(\widehat E_t\neq E_t)\le\sum_{i}\pi_iP^{(i)}(m)+\frac{\bar s}{2}\sum_{k=1}^{m-1}\bigl(\psi_1(k)+\psi_2(k)\bigr)=A(m).
\]

\emph{Lower bound.} A stationary two-state Markov chain is reversible, since \(\pi_1Q_{12}=\pi_2Q_{21}\). The observations are independent given the regime path, with laws that depend only on the current regime, so the window \(((E_s,Y_s))_{t-m<s\le t}\) has the same law as its time reversal. The reversal preserves \(S\), and hence the decision, and exchanges \(E_{t-m+1}\) and \(E_t\). Therefore
\[
\PP\bigl(\widehat E_t\neq E_t,\ E_{t-m+1}\neq E_t\bigr)=\PP\bigl(\widehat E_t\neq E_{t-m+1},\ E_{t-m+1}\neq E_t\bigr).
\]
On \(\{E_{t-m+1}\neq E_t\}\) exactly one of the two events occurs, so both probabilities equal \(\tfrac12\PP(E_{t-m+1}\neq E_t)=\pi_1(Q^{m-1})_{12}=\pi_1\pi_2(1-\lambda^{m-1})\), which is \eqref{eq:reversal}. The window lies entirely in regime~\(i\) with probability \(\pi_iQ_{ii}^{m-1}\), and its error is then \(P^{(i)}(m)\), which gives the first inequality. For the second, the events \(\mathcal A_{k,i}:=\{E_t=i,\ \tau_t=k,\ E_{t-m+1}=\cdots=E_{t-k}\}\), \(1\le k\le m-1\), on which the window contains exactly one switch, are disjoint subsets of \(\{E_{t-m+1}\neq E_t\}\), and with \(j\neq i\)
\[
\PP(\mathcal A_{k,i})=\pi_jQ_{jj}^{m-k-1}Q_{ji}Q_{ii}^{k-1}\ge\frac{\bar s}{2}\,p_{\min}^{m-2}.
\]
Hence \(\pi_1\pi_2(1-\lambda^{m-1})\ge\tfrac12\sum_{k,i}\PP(\mathcal A_{k,i})\ge\tfrac{\bar s}{2}(m-1)p_{\min}^{m-2}\), and \(\pi_iQ_{ii}^{m-1}\ge\pi_ip_{\min}^{m-1}\).

\emph{Upper bound under Markov switching.} Let \(N\) be the number of switches in the window, so that \(E_{t-m+1}\neq E_t\) if and only if \(N\) is odd, and \(\EE N=\bar s(m-1)\). Windows without switches contribute \(\sum_i\pi_iQ_{ii}^{\,m-1}P^{(i)}(m)\le\bar P(m)\). An error together with an odd number of switches has probability \(\tfrac12\PP(N\text{ odd})\) by \eqref{eq:reversal}, and an error together with a positive even number has probability at most \(\PP(N\ge2,\ N\text{ even})\). Since \(N\ge1\) when \(N\) is odd and \(N\ge2\) when \(N\) is positive and even,
\[
\PP(\widehat E_t\neq E_t)\le\bar P(m)+\tfrac12\PP(N\text{ odd})+\PP(N\ge2,\ N\text{ even})\le\bar P(m)+\tfrac12\EE N=A(m).
\]
This uses neither monotonicity nor the pairing.
\end{proof}

The law needs neither symmetric switching nor mirror-image noise. Under Markov switching it needs no ordering of the emission laws either, and part~(i) needs no Markov property and holds for regimes whose sojourns are not geometric, such as bull and bear markets whose ending probabilities depend on their age \cite{LundeTimmermann2004}, which a hidden Markov filter misspecifies (Experiment~I). Gaussian residuals with a common variance and \(\pm1\) observations, for which \(P^{(i)}(m)\) is a binomial tail, satisfy the assumptions. The factor \(p_{\min}^{m-1}\ge1-(1-p_{\min})(m-1)\) is close to one near the optimum, because optimal windows are much shorter than a typical sojourn, and tends to one there as switches become rare, so \(A\) is asymptotically exact at the optimum.

\begin{figure}[htbp]
\centering
\begin{tikzpicture}
\begin{axis}[
    width=0.8\textwidth,
    height=6cm,
    xlabel={Window length \(m\)},
    ylabel={Error probability},
    xmin=1, xmax=60, ymin=0, ymax=0.62,
    axis lines*=left,
    ymajorgrids,
    grid style={line width=0.3pt, draw=black!8},
    axis line style={line width=0.4pt},
    tick align=outside,
    tick style={line width=0.4pt, black},
    tick label style={font=\footnotesize},
    label style={font=\small},
    legend style={draw=none, fill=white, fill opacity=0.9, text opacity=1, font=\footnotesize, row sep=1pt, cells={anchor=west},
                  at={(0.56,0.97)}, anchor=north},
]
\input{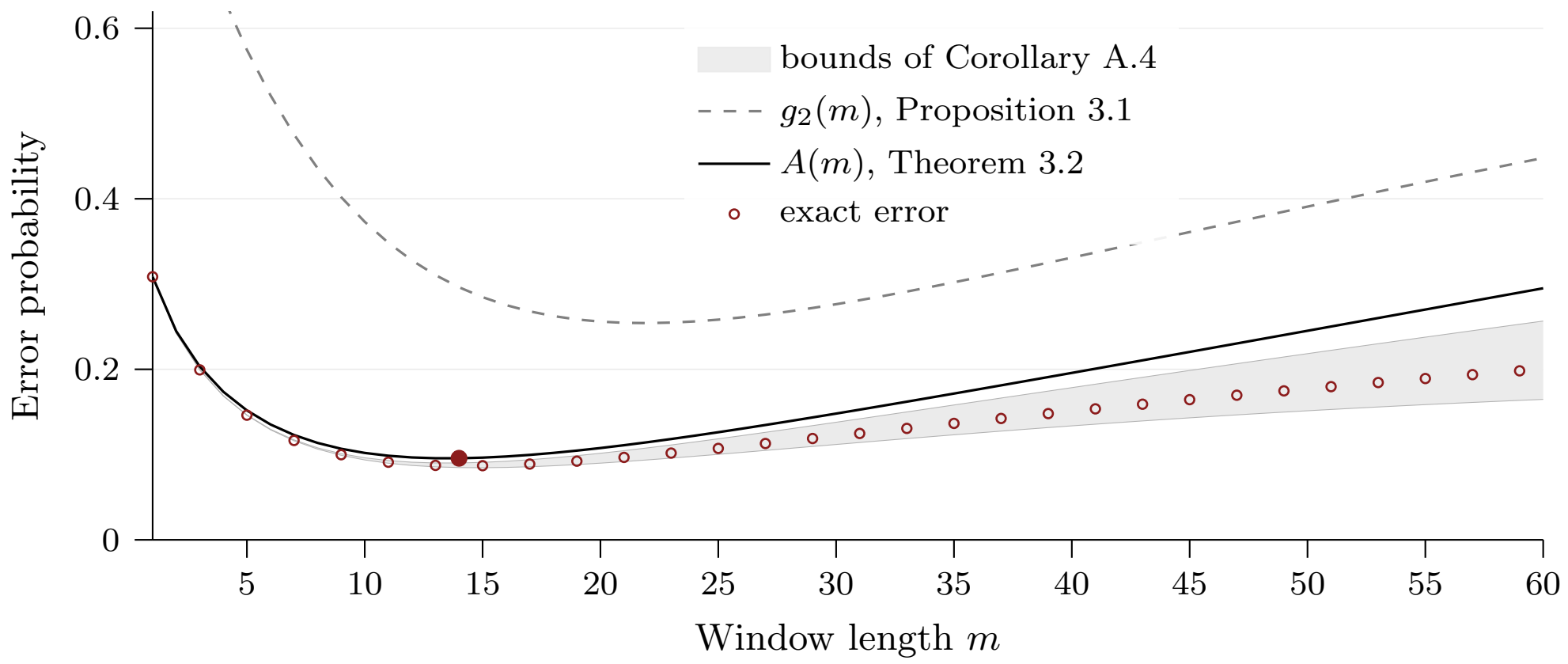}
\end{axis}
\end{tikzpicture}
\caption{Gaussian residuals, \(\Delta=\sigma=1\), \(p=\G{gauss-ill-p}\) (\(\kappa=0.01\)). The law \(A(m)\) is minimized at \(m_A=\G{gauss-ill-mA}\) (filled marker), the optimal window, and \(g_2\) at \(m=\G{gauss-ill-g2arg}\).}
\label{fig:half-window}
\end{figure}

The proof of~(ii) uses only that reversing the window leaves the decision unchanged, so~\eqref{eq:reversal} holds for every such decoder, for example any function of the empirical distribution of the window. Given that the oldest and current regimes differ, such a statistic carries no information about which of the two is current. Part~(ii) does not use the stochastic order, and part~(i) uses it only to show that \(\psi_1\) and \(\psi_2\) are nonincreasing. Monotonicity also holds for Gaussian residuals with unequal standard deviations \(\sigma_1\neq\sigma_2\), which are not ordered (Appendix~\ref{app:proofs-local}), so the law holds with \(P^{(i)}(m)=\Phibar(\sqrt m\,\Delta/(2\sigma_i))\). Away from the optimum the two sides of the law separate. The bound \(A\) treats a window with several switches as if its last switch were the only one, while the lower bound of~(ii) tends to \(\pi_1\pi_2\) for long windows. For the symmetric chain, Corollary~\ref{cor:symmetric-bounds} gives bounds that remain informative at every window length (shaded band in Figure~\ref{fig:half-window}).

\subsection{The optimal window}

\begin{corollary}[Window rule and sharp asymptotics]
\label{cor:window-rule}
Under the assumptions of Theorem~\ref{thm:half-window}(ii), let \(P^\star:=\inf_{m\ge1}\PP(\widehat E_t\neq E_t)\), and let \(m_A\) minimize \(A\), which exists because \(A(m)\to\infty\) when \(\bar s>0\).
\begin{enumerate}[label=(\roman*)]
\item The window \(m_A\) is near-optimal: if \(p_{\min}>0\), its error is at most \(p_{\min}^{-(m-1)}\) times the error of the window of length \(m\), for every \(m\ge1\). In particular, if the infimum \(P^\star\) is attained at \(m^\star\), the error of \(m_A\) is at most \(P^\star/p_{\min}^{\,m^\star-1}\).
\item For Gaussian residuals with common variance \(\sigma^2\), the minimizer over real \(m>0\) of \(\Phibar(\sqrt m\,\Delta/(2\sigma))+\bar s\,m/2\) is
\begin{equation}
m_\kappa=\frac{4\sigma^2}{\Delta^2}\,z_\kappa^2,
\qquad
\frac{\varphi(z_\kappa)}{z_\kappa}=4\kappa,
\label{eq:window-rule}
\end{equation}
with \(\varphi\) the standard normal density and \(\kappa=\bar s\sigma^2/\Delta^2\).
\item Let
\[
I_1:=\sup_{\lambda\ge0}\bigl\{\lambda\theta-\log\EE_{F_1}e^{\lambda Y}\bigr\},
\qquad
I_2:=\sup_{\lambda\ge0}\bigl\{-\lambda\theta-\log\EE_{F_2}e^{-\lambda Y}\bigr\}
\]
be the Chernoff exponents of \(P^{(1)}\) and \(P^{(2)}\), and \(I:=\min\{I_1,I_2\}\). Assume that \(F_1\) and \(F_2\) have finite exponential moments near zero and that \(\theta\) lies in the interior of the convex hull of the support of each, so that \(0<I<\infty\). Let the switching probabilities tend to zero with a fixed ratio and fixed emission laws. Then, for small enough switching probabilities, the infimum \(P^\star\) is attained at some \(m^\star\), and
\[
P^\star=\frac{\bar s\log(1/\bar s)}{2I}\,(1+o(1)),
\qquad
m^\star,\ m_A=\frac{\log(1/\bar s)}{I}\,(1+o(1)).
\]
For Gaussian residuals with common variance, \(I=\Delta^2/(8\sigma^2)\), so \(P^\star=4\kappa\log(1/\kappa)(1+o(1))\) and \(m^\star,m_A,m_\kappa=(8\sigma^2/\Delta^2)\log(1/\kappa)(1+o(1))\).
\end{enumerate}
\end{corollary}

Measured in resolution times, the rule \eqref{eq:window-rule} depends only on \(\kappa\), which explains the collapse of the optimal windows and errors in Figure~\ref{fig:expII-collapse}. Part~(iii) is proved in Appendix~\ref{app:proofs-local}. Windows much longer than a sojourn are not optimal when switches are rare, since by Theorem~\ref{thm:half-window}(ii) their error is at least \(\pi_1\pi_2(1-\lambda^{m-1})\). Windows of length \(O(1/\bar s)\) are compared by bounding their errors from below at once (Lemma~\ref{lem:three-ranges}).

\emph{Choosing the window from data.}\label{sec:plug-in} The rule needs \(p\), \(\Delta\) and \(\sigma\), which in practice are estimated from a training trajectory with known regimes, as in \cite{PambukyanSaudagarKumar2026}. Assume that the training trajectory visits both regimes and contains at least one switch. Estimating \(p\) by \(\widehat p:=1-N/(T-1)\), where \(N\) is the number of switches in the training trajectory, and \(\mu_i\), \(\sigma\) by regime-wise sample means and the pooled residual standard deviation, gives the plug-in rule
\begin{equation}
\widehat m:=\arg\min_m\Bigl\{\bar\Phi\bigl(\sqrt m\,(\widehat\mu_2-\widehat\mu_1)/(2\widehat\sigma)\bigr)+(1-\widehat p)\,\frac{m-1}{2}\Bigr\}.
\label{eq:plug-in}
\end{equation}
For the symmetric chain \(N\) is binomial, so the relative error of \(1-\widehat p\) is of order \(N^{-1/2}\) as the expected number of switches \(T\bar s\) grows, so its accuracy is set by the number of switches rather than by the trajectory length. The error curve is flat near its minimum, so a misestimated \(\kappa\) costs little. In Experiment~I, using \(\kappa/2\) or \(2\kappa\) in the rule \eqref{eq:window-rule} raises the error by at most \(\G{gauss-robust-half}\%\) and \(\G{gauss-robust-double}\%\) when \(\kappa\le0.03\). Underestimating \(\kappa\) lengthens the window and is the safer error.

\subsection{Comparison with the Bayes filter}

As the switching probabilities tend to zero with a fixed ratio, the error of the Bayes filter is
\begin{equation}
P_{\mathrm{Bayes}}=\frac{\bar s}{2}\Bigl(\frac{1}{\mathrm{KL}(F_2\,\|\,F_1)}+\frac{1}{\mathrm{KL}(F_1\,\|\,F_2)}\Bigr)\log(1/\bar s)\,(1+o(1)),
\label{eq:kz}
\end{equation}
by the two-state case of Theorem~1 of Khasminskii and Zeitouni \cite{KhasminskiiZeitouni1996}, for mutually absolutely continuous emission laws whose log-likelihood ratio has a finite moment of order larger than two. Comparing with Corollary~\ref{cor:window-rule}(iii), the optimal uniform window loses the factor
\begin{equation}
\frac{P^\star}{P_{\mathrm{Bayes}}}\;\longrightarrow\;\frac{1}{I\,\bigl(1/\mathrm{KL}(F_2\,\|\,F_1)+1/\mathrm{KL}(F_1\,\|\,F_2)\bigr)},
\label{eq:price}
\end{equation}
which is at least one because the filter is optimal, and equals \(\mathrm{KL}/(2I)\) for mirror-image laws. For Gaussian residuals with common variance \(\mathrm{KL}=\Delta^2/(2\sigma^2)\) and \(I=\Delta^2/(8\sigma^2)\), so \(P_{\mathrm{Bayes}}=2\kappa\log(1/\kappa)(1+o(1))\) and uniform windows lose a factor two asymptotically. For \(\pm1\) observations with mirror-image means \(\mp\Delta/2\) and \(q:=(1-\Delta/2)/2\) the factor is
\[
\frac{(1-2q)\log\frac{1-q}{q}}{\log\frac{1}{4q(1-q)}},
\]
which is close to two for small \(\Delta\) (\(1.94\) at \(\Delta=0.8\)) and decreases to one as \(\Delta\to2\). For Gaussian residuals the factor two comes from the shape of the window, as the next proposition shows.

\begin{proposition}[Window shape]
\label{prop:shape}
Let \(g:[0,\infty)\to[0,\infty)\) be nonincreasing with \(\int g=1\). Let \(u_g\) satisfy \(\int_0^{u_g}g=1/2\), and set
\[
c_g:=u_g\int_0^\infty g^2 .
\]
For \(L\ge1\), let \(P_L:=\limsup_{t\to\infty}\PP(\widehat E_t\neq E_t)\) be the error of the decoder with trailing weights \(w_{t-j}\propto g(j/L)\), \(0\le j<t\), on the stationary symmetric chain under Assumption~\ref{ass:residual-model}. As \(p\to1\) with \(\Delta/\sigma\) fixed,
\[
\inf_{L\ge1}P_L\le 8c_g\,\kappa\log(1/\kappa)\,(1+o(1)),
\]
and for independent \(N(0,\sigma^2)\) residuals the infimum over \(1\le L\le C/(1-p)\), for any fixed \(C\), equals \(8c_g\,\kappa\log(1/\kappa)\,(1+o(1))\). Moreover \(c_g\ge1/4\) for every \(g\), uniform windows (\(g=\mathbf 1_{[0,1)}\)) have \(c_g=1/2\), and exponential weights (\(g(u)=e^{-u}\)) have \(c_g=\tfrac12\log2\).
\end{proposition}

The constant \(8c_g\) is \(4\) for uniform windows and \(4\log2\) for exponential weights, against \(2\) for the Bayes filter \eqref{eq:kz}. Uniform windows therefore lose a factor \(2\). Exponential weights lose at most a factor \(2\log2\approx1.39\), and exactly this factor among scales \(L\) of order at most \(1/(1-p)\), the mean sojourn time. Longer scales are not excluded by the proof, although we expect them to be worse. A switch is detected once the recent part of the window carries half of the weight, which takes \(u_gL\) steps, while the noise is controlled by \(\|w\|^2\approx v_g/L\) with \(v_g:=\int g^2\). Under uniform weights the older half of the window adds noise without helping to detect the switch, and exponential weights reduce this loss. The bound \(c_g\ge1/4\), from the Cauchy--Schwarz inequality, matches the Bayes rate but is attained by no fixed shape, since equality would require \(g\) to vanish beyond \(u_g\). The proposition holds for each fixed \(g\), and its \(o(1)\) term depends on \(g\). Approaching \(1/4\) requires a shape that changes with \(\kappa\). At practical \(\kappa\) the differences are smaller than the limits (Experiment~I).

Allowing the shape to depend on \(\kappa\) attains the Bayes rate in the limit. A recent block detects a switch once it carries half of the weight, and an older block spread thinly over a longer stretch adds little noise. The length of the older block balances that noise against the slowly growing margin after a switch.

\begin{proposition}[Two-block windows]
\label{prop:two-block}
For the stationary symmetric chain under Assumption~\ref{ass:residual-model}, let \(\ell:=\log(1/\kappa)\), \(r:=(4\ell/\pi)^{1/4}\) with \(\pi=3.14\ldots\), \(a:=\lceil2(1+1/r)\,\ell\,\sigma^2/\Delta^2\rceil\) and \(b:=\lceil ra\rceil\), and give the weight \(1/(2a)\) to each of the \(a\) most recent observations and \(1/(2b)\) to each of the \(b\) before them. For \(t\ge a+b\), as \(p\to1\) with \(\Delta/\sigma\) fixed,
\[
\PP(\widehat E_t\neq E_t)\le2\kappa\log(1/\kappa)\Bigl(1+2\Bigl(\frac{\pi}{4\log(1/\kappa)}\Bigr)^{1/4}+O\bigl((\log(1/\kappa))^{-1/2}\bigr)\Bigr).
\]
For independent \(N(0,\sigma^2)\) residuals this is the Bayes rate \eqref{eq:kz} to first order.
\end{proposition}

The bound exceeds the leading term \(2\kappa\log(1/\kappa)\) of the Bayes error by a relative amount of order \((\log(1/\kappa))^{-1/4}\). Since \eqref{eq:kz} gives the Bayes error only to first order, this is not a rate relative to the Bayes error itself, and it decays slowly, so the proposition describes the limit rather than practical switching rates. The proof in Appendix~\ref{app:proofs-local} also gives an explicit bound for every \(\kappa\in(0,1)\).

When the decision can wait, the window can be centred on \(t\). A switch then causes an error only if it falls near the centre, and the switching cost drops from order \(m\) to order \(\sqrt m\).

\begin{proposition}[Centred windows]
\label{prop:centred}
For the stationary symmetric chain under Assumption~\ref{ass:residual-model} and the centred uniform window \(\{t-h,\ldots,t+h\}\), \(m=2h+1\), \(t>h\),
\[
\PP(\widehat E_t\neq E_t)\le e^{-m\Delta^2/(8\sigma^2)}+2(1-p)\left(1+\frac{\sigma}{\Delta}\sqrt{\frac{\pi m}2}\right)+\frac{(1-p)^2(m-1)^2}{2}.
\]
With \(m\) the smallest odd integer at least \((8\sigma^2/\Delta^2)\log(1/\kappa)\), the right-hand side is \(4\sqrt\pi\,\kappa\sqrt{\log(1/\kappa)}\,(1+o(1))\) as \(p\to1\).
\end{proposition}

Propositions~\ref{prop:shape}--\ref{prop:centred} are proved in Appendix~\ref{app:proofs-local}.

\begin{remark}[Regret]
\label{rem:regret}
Window decoders are often used to act, as in moving-average trading rules. Suppose an agent chooses one of two actions from \(Y_1,\ldots,Y_{t-1}\), action~\(i\) is optimal in regime~\(i\), and the other action costs one unit more, as when holding an asset long or short whose price is an MMRW with drifts \(\mp\Delta/2\). The action defines a decoder \(\widehat E_{t-1}\). For the symmetric chain with \(p>1/2\), \(E_t=E_{t-1}\) with probability \(p\) independently of the past observations, so the expected per-step regret against an agent who knows \(E_t\) is
\[
\PP(\widehat E_{t-1}\neq E_t)=(1-p)+(2p-1)\,\PP\bigl(\widehat E_{t-1}\neq E_{t-1}\bigr).
\]
The avoidable part of the regret is therefore proportional to the decoding error, and the factors of this section apply to it. For Gaussian residuals with common variance a moving-average rule has asymptotically twice the avoidable regret of an agent who acts on the Bayes filter, and an exponentially weighted rule at most \(2\log2\) times as much. For an asymmetric chain the regret is affine in the two kinds of error separately rather than in their sum.
\end{remark}

\section{More Than Two Regimes: Switch Costs}\label{sec:k-regimes}

With \(K\) regimes, a window that contains a switch can be assigned to a third regime. As the switch moves through the window, the mean of the window average travels from the old regime mean to the new one and crosses the decision regions of the regimes in between. Theorem~\ref{thm:switch-cost} extends the half-window law by assigning a cost to each switch.

\subsection{The switch-cost law}

As in Theorem~\ref{thm:half-window}, the observations are independent given the regime path, with \(Y_s\sim F_{E_s}\), and the decoder depends on the window only through the sum of its observations, with ties broken at random, as \eqref{eq:decoder} does. For \(i\neq j\) and \(1\le k\le m-1\), let \(e_{ji}(k)\) be its error when regime~\(i\) is current, the \(k\) most recent observations come from regime~\(i\) and the other \(m-k\) from regime~\(j\), and let \(\chi_{ij}(k)\) be the probability that it returns neither \(i\) nor \(j\) for this composition of the window. The \emph{cost} of a switch from \(j\) to \(i\) and the rate of such switches are
\[
c_{ji}(m):=\frac1{m-1}\sum_{k=1}^{m-1}e_{ji}(k),
\qquad
\bar s_{ji}:=\PP(E_{s-1}=j,\ E_s=i),
\]
so that \(\bar s=\sum_{i\neq j}\bar s_{ji}\), and \(\bar s_{ji}=\pi_jQ_{ji}\) for a Markov chain. With \(P^{(i)}(m)\) and \(\bar P(m)\) as in Theorem~\ref{thm:half-window}, let
\begin{equation}
A_K(m):=\bar P(m)+(m-1)\sum_{i\neq j}\bar s_{ji}\,c_{ji}(m).
\label{eq:AK}
\end{equation}
For \(m=1\) the sums are empty, and we set \(c_{ji}(1):=0\) and \(\chi_{ij}(1):=0\), so that \(A_K(1)=\bar P(1)\).

\begin{theorem}[Switch-cost law]
\label{thm:switch-cost}
Let \((E_t)\) be a stationary process on \(K\) regimes, and let \(N_{t,m}\) be the number of switches in the trailing window of length \(m\le t\).
\begin{enumerate}[label=(\roman*)]
\item \(\PP(\widehat E_t\neq E_t)\le A_K(m)+\PP(N_{t,m}\ge2)\).
\item If \((E_t)\) is a Markov chain, then \(\PP(\widehat E_t\neq E_t)\ge p_{\min}^{\,m-1}A_K(m)\) and \(\PP(N_{t,m}\ge2)\le\binom{m-1}{2}(1-p_{\min})\,\bar s\).
\item For \(m\ge2\) and \(i\neq j\), \(c_{ji}(m)+c_{ij}(m)=1+\chi_{ij}(m)\), where \(\chi_{ij}(m):=\frac1{m-1}\sum_{k=1}^{m-1}\chi_{ij}(k)\). Consequently, if \(\bar s_{ij}=\bar s_{ji}\) for every pair, as for a reversible chain,
\[
A_K(m)=\bar P(m)+\bar s\,\frac{m-1}{2}+(m-1)\sum_{i<j}\bar s_{ij}\,\chi_{ij}(m).
\]
\end{enumerate}
\end{theorem}

\begin{proof}
(i) Split according to the number of switches in the window. If there is none and \(E_t=i\), the error is \(P^{(i)}(m)\), and this event has probability at most \(\pi_i\). If there is exactly one, regime~\(j\) occupies the first \(m-k\) positions of the window and regime~\(i\) the last \(k\), for some \(j\neq i\) and \(1\le k\le m-1\). The error is then \(e_{ji}(k)\), and the event requires a switch from \(j\) to \(i\) at time \(t-k+1\), which has probability \(\bar s_{ji}\) by stationarity. Bounding the error by one when there are two or more switches,
\[
\PP(\widehat E_t\neq E_t)\le\sum_i\pi_iP^{(i)}(m)+\sum_{i\neq j}\bar s_{ji}\sum_{k=1}^{m-1}e_{ji}(k)+\PP(N_{t,m}\ge2)=A_K(m)+\PP(N_{t,m}\ge2).
\]

(ii) For a Markov chain, the two kinds of windows of~(i) have probabilities
\[
\pi_iQ_{ii}^{m-1}\ge\pi_i\,p_{\min}^{m-1}
\qquad\text{and}\qquad
\pi_jQ_{jj}^{m-k-1}Q_{ji}Q_{ii}^{k-1}\ge\bar s_{ji}\,p_{\min}^{m-2},
\]
and dropping the windows with two or more switches gives the lower bound. Given the path up to a transition, the transition is a switch with probability at most \(1-p_{\min}\). Two given transitions of the window are therefore both switches with probability at most \(\bar s(1-p_{\min})\), and there are \(\binom{m-1}{2}\) pairs.

(iii) A window with \(k\) observations from regime~\(i\) and \(m-k\) from regime~\(j\) gives the same law of the sum, and hence of the decoder's output \(D\), whether \(i\) is current after a switch from \(j\), or \(j\) is current after a switch from \(i\) at the mirrored position. The decoder errs in the first case if \(D\neq i\) and in the second if \(D\neq j\), so
\[
e_{ji}(k)+e_{ij}(m-k)=\PP(D\neq i)+\PP(D\neq j)=1+\chi_{ij}(k).
\]
Summing over \(k\) gives the identity. If \(\bar s_{ij}=\bar s_{ji}\), summing over pairs and using \(\sum_{i<j}\bar s_{ij}=\bar s/2\) gives the last formula.
\end{proof}

For two regimes, \(\chi_{12}\equiv0\) and \(\bar s_{12}=\bar s_{21}=\bar s/2\) for every stationary process, so \(A_2=A\) and the half-window law is the case of a single pair of regimes. Theorem~\ref{thm:half-window}(i) also removes the term \(\PP(N_{t,m}\ge2)\), by monotonicity. With more regimes this argument fails, because the observations before the last switch may come from a regime farther away than the previous one. The term is of second order in \(\bar s m\) when the switching probabilities are of comparable size.

\subsection{Switch costs and the regime means}

Let \(\mu_1<\cdots<\mu_K\), let \(R_i\) be the decision region of regime~\(i\) under \eqref{eq:decoder}, the interval around \(\mu_i\) bounded by the midpoints between \(\mu_i\) and its neighbours, and for \(i\neq j\) let \(h_{ij}\) be the distance from \(\mu_i\) to the end point of \(R_i\) facing \(\mu_j\).

\begin{proposition}[Geometry of a switch]
\label{prop:switch-geometry}
Let the emission laws have means \(\mu_1<\cdots<\mu_K\) and finite variances, and consider the decoder \eqref{eq:decoder}. As \(m\to\infty\), for every pair \(i\neq j\),
\[
c_{ji}(m)\to c_{ji}:=1-\frac{h_{ij}}{|\mu_i-\mu_j|},
\qquad
\chi_{ij}(m)\to1-\frac{h_{ij}+h_{ji}}{|\mu_i-\mu_j|}.
\]
Hence \(c_{ji}\ge1/2\), with equality if and only if no mean lies between \(\mu_i\) and \(\mu_j\). For Gaussian residuals with common variance,
\[
\sum_{k=1}^{m-1}e_{ji}(k)=m\,c_{ji}-\tfrac12+\epsilon_m,
\]
where \(\epsilon_m\to0\) faster than any power of \(m\).
\end{proposition}

The proof is in Appendix~\ref{app:proofs-k}. A switch between neighbouring regimes costs half a window, as with two regimes, and a switch over intermediate regimes costs more, since for a fraction \(1-(h_{ij}+h_{ji})/|\mu_i-\mu_j|\) of its positions the window average lies in their decision regions, where it is wrong for both the old and the new regime. With equally spaced means, a jump across \(d\) gaps costs \(1-1/(2d)\) of a window. Weighting the costs by the switch rates gives the \emph{mean cost of a switch},
\[
\bar c:=\frac1{\bar s}\sum_{i\neq j}\bar s_{ji}\,c_{ji}\in\bigl[\tfrac12,1\bigr),
\]
and for Gaussian residuals the proposition gives \(A_K(m)=\bar P(m)+\bar c\,\bar s\,m-\bar s/2+\bar s\,\epsilon'_m\), with \(\epsilon'_m\to0\) faster than any power of \(m\). Up to the remainder \(\bar s\,\epsilon'_m\), which can depend on the individual rates \(\bar s_{ji}\), the law depends on the switching only through \(\pi\), \(\bar s\) and \(\bar c\).

The cost depends on the direction of a switch when the gaps differ, since \(h_{ij}\) is measured on the side of arrival, so the error of a window can change when a non-reversible environment is reversed in time (Experiment~II). The same argument applies to vector observations and the nearest-mean decoder, whose regions are Voronoi cells. Then \(c_{ji}\) is the fraction of the segment from \(\mu_j\) to \(\mu_i\) that lies outside the cell of \(\mu_i\), which equals \(1/2\) exactly when no other mean lies in the open ball with diameter \([\mu_j,\mu_i]\).

\subsection{The optimal window and comparison with the Bayes filter}

\begin{corollary}[Window rule for \(K\) regimes]
\label{cor:k-window}
Let \((E_t)\) be a stationary Markov chain, let the observations be as in Proposition~\ref{prop:switch-geometry}, let \(P^\star\) and \(m^\star\) be as in Corollary~\ref{cor:window-rule}, and let \(m_{A_K}\) minimize \(A_K\).
\begin{enumerate}[label=(\roman*)]
\item If \(p_{\min}>0\) and the infimum \(P^\star\) is attained at \(m^\star\), the error of the window \(m_{A_K}\) is at most \(P^\star/p_{\min}^{\,m^\star-1}+\binom{m_{A_K}-1}{2}(1-p_{\min})\,\bar s\).
\item Let \(Q=\mathrm{Id}+\varepsilon\Gamma\) for a fixed irreducible generator \(\Gamma\), and let the emission laws be fixed, with finite exponential moments near zero. Let \(I\) be the smallest Chernoff exponent of the single-regime errors, over the regimes and the end points of their decision regions, and assume that every end point of \(R_i\) lies in the interior of the convex hull of the support of \(F_i\), so that \(0<I<\infty\). As \(\varepsilon\to0\), the infimum \(P^\star\) is eventually attained, and
\[
P^\star=\frac{\bar c\,\bar s\log(1/\bar s)}{I}\,(1+o(1)),
\qquad
m^\star,\ m_{A_K}=\frac{\log(1/\bar s)}{I}\,(1+o(1)).
\]
For Gaussian residuals with common variance, \(I=\Delta^2/(8\sigma^2)\), with \(\Delta\) the smallest gap between the means as in \eqref{eq:Delta}, and \(P^\star=8\bar c\,\kappa\log(1/\kappa)\,(1+o(1))\).
\end{enumerate}
\end{corollary}

For two regimes \(\bar c=1/2\), and the corollary reduces to Corollary~\ref{cor:window-rule}(iii). Its proof is in Appendix~\ref{app:proofs-k}. The stationary law \(\pi\) and the mean cost \(\bar c\) do not depend on \(\varepsilon\), \(\bar s\) and \(1-p_{\min}\) are both proportional to \(\varepsilon\), and the costs \(c_{ji}(m)\) depend only on \(m\) and the emission laws, so their convergence in Proposition~\ref{prop:switch-geometry} is uniform in \(\varepsilon\). This uniformity allows the proof to compare all windows of length \(O(1/\varepsilon)\) at once (Lemma~\ref{lem:three-ranges}), and Lemma~\ref{lem:long-windows} excludes longer windows.

For the Bayes filter, Theorem~1 of Khasminskii and Zeitouni \cite{KhasminskiiZeitouni1996} covers any irreducible aperiodic chain with transition probabilities \(\varepsilon\Gamma_{ij}\), \(i\neq j\), including chains that forbid some transitions, under the conditions on the emission laws stated after \eqref{eq:kz}, which Gaussian residuals satisfy:
\begin{equation}
P_{\mathrm{Bayes}}=\sum_{i\neq j}\frac{\bar s_{ji}}{\mathrm{KL}(F_i\,\|\,F_j)}\,\log(1/\bar s)\,(1+o(1)).
\label{eq:kz-K}
\end{equation}
Thus a switch from \(j\) to \(i\) costs the filter \(\log(1/\bar s)/\mathrm{KL}(F_i\,\|\,F_j)\) erroneous steps, against \(c_{ji}\log(1/\bar s)/I\) for the optimal window. For Gaussian residuals with common variance and equally spaced consecutive means with gap \(\Delta\), a switch across \(d\) gaps contributes
\[
\Bigl(1-\frac1{2d}\Bigr)\frac{8\sigma^2}{\Delta^2}\,\log(1/\bar s)
\qquad\text{and}\qquad
\frac{2\sigma^2}{d^2\Delta^2}\,\log(1/\bar s)
\]
erroneous steps to the window and to the filter, a ratio of \(2d(2d-1)\), which is two between neighbours, as for two regimes, and twelve for a jump over one regime. The overall ratio of \(P^\star\) to \eqref{eq:kz-K} is the average of these per-switch ratios weighted by the filter's costs, so it depends on the mix of switches. After a large jump the window average passes through the decision regions of the intermediate regimes, whereas the filter detects the jump quickly.

\section{Application: The Reflected Walk}\label{sec:reflected}

In the reflected MMRW itself, the conditional drift of an increment depends on the position, because near a wall a step towards the wall is blocked. We consider nearest-neighbour kernels with
\[
P^{(e)}_{x,x+1}=u_e,\qquad P^{(e)}_{x,x-1}=d_e\qquad(0<x<M),
\]
the remaining mass on staying put, and two kinds of walls. At \emph{holding} walls a blocked step leaves the walk in place: \(P^{(e)}_{0,1}=u_e\) and \(P^{(e)}_{M,M-1}=d_e\), with the remaining mass on staying put. At \emph{bouncing} walls the walk always moves away: \(P^{(e)}_{0,1}=P^{(e)}_{M,M-1}=1\). The regimes are indexed so that \(E_t\) governs the step from \(X_{t-1}\) to \(X_t\), and the environment is autonomous: \(\PP(E_{t+1}=j\mid E_{0:t},X_{0:t})=Q_{E_tj}\). Write
\[
d(x,e):=\EE[Y_t\mid X_{t-1}=x,\ E_t=e],
\qquad
\mathcal G_s:=\sigma(X_0,\ldots,X_{s-1},E_1,\ldots,E_s),
\qquad
\eta_s:=Y_s-d(X_{s-1},E_s),
\]
so that \(d(x,e)=u_e-d_e=:\mu_e\) in the interior and \(\EE[\eta_s\mid\mathcal G_s]=0\), where \(\eta_s\) is \(\mathcal G_{s+1}\)-measurable.

\subsection{Raw increments}

Visits to a wall can reduce the accuracy of the decoder \eqref{eq:decoder} applied to the raw increments.

\begin{proposition}[Raw increments]
\label{prop:reflected}
Let \(H_{t,m}\) be the event that some step in the window starts at \(0\) or \(M\). For the decoder \eqref{eq:decoder} applied to \(Y_s=X_s-X_{s-1}\) with the interior drifts \(\mu_e\), for every \(t\ge m\ge1\),
\begin{equation}
\PP(\widehat E_t\neq E_t)
\le
2\exp\!\left(-\frac{m\Delta^2}{8}\right)
+\PP(\mathrm{NoSwitch}_{t,m}^c)+\PP(H_{t,m}).
\label{eq:reflected-bound}
\end{equation}
\end{proposition}

The proof is that of Proposition~\ref{thm:recoverability}, with Hoeffding's lemma in place of the sub-Gaussian assumption (Appendix~\ref{app:proofs-reflected}). The boundary term cannot be avoided with raw increments. Since \(\bar Y_{t,m}=(X_t-X_{t-m})/m\) and \(0\le X_s\le M\), \(|\bar Y_{t,m}|\le M/m\) whatever the regimes, so long windows of raw increments cannot separate the drifts of a small lattice. In Experiment~III the bound \eqref{eq:reflected-bound} never falls below \(\G{refl-minbound-lo}\).

\subsection{Position-aware decoding and the information clock}

We now restrict to \(K=2\). The decoder can reduce this loss by using the position of the walk. Let \(\Delta_x:=d(x,2)-d(x,1)\), assumed nonnegative for every \(x\) and positive in the interior, let \(v_x\) be a sub-Gaussian variance proxy for \(Y_t\) given \(X_{t-1}=x\) (a quarter of the squared length of its range), and define the per-step signal-to-noise ratio \(\mathrm{snr}_x:=\Delta_x^2/v_x\), with \(\mathrm{snr}_x:=0\) when \(v_x=0\). The position-aware statistic weights and centres each increment according to its starting position,
\begin{equation}
V_t:=\sum_{s=t-m+1}^{t}\omega_{X_{s-1}}\Bigl(Y_s-\frac{d(X_{s-1},1)+d(X_{s-1},2)}{2}\Bigr),
\qquad
\omega_x:=\frac{\Delta_x}{v_x}\ (\omega_x:=0\text{ if }v_x=0),
\label{eq:position-aware}
\end{equation}
and the decoder chooses regime~2 if \(V_t>0\) and regime~1 if \(V_t<0\), with ties broken by a fair coin. Steps that cannot distinguish the regimes receive no weight. Unlike the decoder on raw increments, this decoder requires the position-dependent drifts \(d(x,e)\) and noise proxies \(v_x\), but still neither the transition matrix nor the full law of the increments. The information the window carries is measured by the \emph{information clock}
\[
\mathcal I_{t,m}:=\sum_{s=t-m+1}^{t}\mathrm{snr}_{X_{s-1}},
\]
which advances by \(\mathrm{snr}_x\) for each step that starts at \(x\).

\begin{proposition}[Position-aware decoding]
\label{prop:position-aware}
For every \(t\ge m\ge1\) and every \(\iota>0\),
\begin{equation}
\PP(\widehat E_t\neq E_t)\le e^{-\iota/8}+\PP(\mathcal I_{t,m}<\iota)+\PP(\mathrm{NoSwitch}_{t,m}^c).
\label{eq:information-clock}
\end{equation}
In particular, if \(\mathrm{snr}_x\ge\mathrm{snr}_{\min}>0\) for every \(x\), then \(\PP(\widehat E_t\neq E_t)\le e^{-m\,\mathrm{snr}_{\min}/8}+\PP(\mathrm{NoSwitch}_{t,m}^c)\).
\end{proposition}

\begin{proof}
Fix the regime path \(E_{t-m+1}=\cdots=E_t=1\), the case of regime~2 being symmetric. On this path each summand of \(V_t\) equals \(\omega_s(\eta_s-\Delta_{X_{s-1}}/2)\), with \(\omega_s:=\omega_{X_{s-1}}\) measurable with respect to \(\mathcal G_s\). By Hoeffding's lemma, with \(\lambda=1/2\),
\[
\EE\bigl[e^{\lambda\omega_s(\eta_s-\Delta_{X_{s-1}}/2)}\bigm|\mathcal G_s\bigr]
\le\exp\Bigl(-\frac{\mathrm{snr}_{X_{s-1}}}{4}+\frac{\mathrm{snr}_{X_{s-1}}}{8}\Bigr)=e^{-\mathrm{snr}_{X_{s-1}}/8}.
\]
Since \(\mathrm{snr}_{X_{s-1}}\) is also \(\mathcal G_s\)-measurable, conditioning successively as in Lemma~\ref{lem:path-mgf} and using the autonomy of the environment shows that
\[
Z_r:=\mathbf 1\{E_{t-m+1:r}=1\}\exp\Bigl(\sum_{s=t-m+1}^{r}\Bigl[\lambda\omega_s\Bigl(\eta_s-\frac{\Delta_{X_{s-1}}}{2}\Bigr)+\frac{\mathrm{snr}_{X_{s-1}}}{8}\Bigr]\Bigr),
\qquad t-m+1\le r\le t,
\]
satisfies \(\EE[Z_r]\le\PP(E_{t-m+1:r}=1)\). On the event \(\{E_{t-m+1:t}=1,\ V_t\ge0,\ \mathcal I_{t,m}\ge\iota\}\) we have \(Z_t\ge e^{\iota/8}\), so this event has probability at most \(e^{-\iota/8}\PP(E_{t-m+1:t}=1)\). Summing over the two regimes and adding \(\PP(\mathcal I_{t,m}<\iota)\) and \(\PP(\mathrm{NoSwitch}_{t,m}^c)\) gives \eqref{eq:information-clock}. For the second claim, take \(\iota=m\,\mathrm{snr}_{\min}\).
\end{proof}

We consider two boundary mechanisms. With holding walls and \(\pm1\) attempted steps, \(u_e=(1+\mu_e)/2=1-d_e\). In the interior \(\Delta_x=\Delta\) and \(v_x=1\). At \(0\) the increment is \(0\) or \(1\), with \(\Delta_0=\Delta/2\) and \(v_0=1/4\), and likewise at \(M\). Hence \(\mathrm{snr}_x=\Delta^2\) everywhere, and the walls change neither the signal-to-noise ratio of a step nor \(\kappa\). For \(\mu_1=-\mu_2\) the summand of \eqref{eq:position-aware} is \(\Delta Y_s\) in the interior and \(2\Delta(Y_s\mp\tfrac12)\) at \(0\) and \(M\), which in each case is \(\Delta\) times the attempted step, so the decoder is the free-walk decoder applied to the attempted steps, in which a blocked step is visible as a zero increment.

With bouncing walls the increment at a wall is the same under both regimes, \(\mathrm{snr}_0=\mathrm{snr}_M=0\), and the clock advances only in the interior, so \(\mathcal I_{t,m}\) is \(\Delta^2\) times the number of interior steps in the window.

\section{Numerical Experiments}\label{sec:experiments}

The experiments test each result against exact computation where one is available.\footnote{The code that generates every figure, table and number in this section is archived as version~1.0.0 \cite{PambukyanNersisyanCode2026} and maintained at \url{https://github.com/arinapambukyan/sliding-window-regime-detection}.} The stationary error of the uniform window is computed exactly by a forward recursion over the regime and the composition of the window, and the boundary probabilities of the reflected walk exactly from the joint chain \((X_t,E_t)\). The stationary error of the Bayes filter for two regimes is approximated numerically. Its predicted log-odds is a Markov chain on a bounded interval, whose kernel we discretize on \(\G{gauss-grid}\) cells. Doubling the number of cells changes the error by at most \(\G{gauss-grid-dev}\%\) in every configuration. Simulations cover the decoders and switching laws without a closed form.

\emph{Monte Carlo error and tuning.} Each simulated error averages the error rates of \(R\) independent stationary trajectories. Its standard error is the standard deviation of these rates divided by \(\sqrt R\), which accounts for the dependence within a trajectory, and ratios carry delta-method standard errors. Standard errors are given in parentheses or bounded once for a group of values. Tuning parameters, the effective length of exponential weights and the half-width of a centred window, are chosen by grid search on separate trajectories, and the errors are estimated on fresh ones. The uniform window is evaluated at the exact optimum or at the rules of Sections~\ref{sec:half-window} and~\ref{sec:k-regimes}, and needs no tuning.

\subsection{Simulation models}

Unless stated otherwise the regime chain is the symmetric two-state chain
\begin{equation}
Q(p)=\begin{pmatrix}p&1-p\\1-p&p\end{pmatrix},
\label{eq:sim-Q}
\end{equation}
started in stationarity, so that \(\bar s=1-p\) and \(\kappa=(1-p)\sigma^2/\Delta^2\). The regime means are \(\mu_1=-\Delta/2\) and \(\mu_2=\Delta/2\), and there are three observation models.
\begin{enumerate}[label=(\alph*)]
\item \emph{\(\pm1\) observations.} \(Y_t\in\{-1,+1\}\), independent given the regimes, with \(\PP(Y_t=+1\mid E_t=i)=(1+\mu_i)/2\). The residuals take values in an interval of length two, so \(\sigma=1\) by Hoeffding's lemma.
\item \emph{Gaussian residuals.} \(Y_t=\mu_{E_t}+\sigma_{E_t}Z_t\) with \(Z_t\) independent standard normal, and \(\sigma_1=\sigma_2=\sigma\) unless stated otherwise.
\item \emph{The reflected walk.} \(X_t=\min\{\max\{X_{t-1}+Z_t,0\},M\}\), with attempted steps \(Z_t\in\{-1,+1\}\) drawn as in~(a). This is the walk of Section~\ref{sec:reflected} with \(u_e=(1+\mu_e)/2\) and holding walls. With bouncing walls the walk moves from \(0\) to \(1\) and from \(M\) to \(M-1\).
\end{enumerate}
Models (a) and (b) satisfy the assumptions of Theorems~\ref{thm:half-window} and~\ref{thm:switch-cost}, and with unequal noise levels the half-window law still holds (Section~\ref{sec:half-window}).

\subsection{Experiment I: two regimes}\label{sec:exp-two}

\emph{Gaussian residuals.} We compute the exact error of the trailing uniform window for every window length in \(\G{gauss-nconfig}\) symmetric configurations, with \(\kappa\) between \(\G{gauss-kappa-lo}\) and \(\G{gauss-kappa-hi}\) and four signal-to-noise ratios \(\Delta/\sigma\), together with the numerically approximated error of the Bayes filter. The same grid is repeated for asymmetric chains with \(Q_{21}=\G{gauss-skew-ratio}\,Q_{12}\) and the same switch rate, and two further configurations have unequal noise levels. At \(\Delta=\sigma=1\) and four values of \(\kappa\), exponential weights, centred windows and the forward--backward smoother are simulated on \(\G{gauss-Rsim}\) trajectories of length \(\G{gauss-Tsim}\), with tuning parameters chosen on \(\G{gauss-Rtune}\) separate trajectories. All standard errors are at most \(\G{gauss-sim-se-max}\).

\begin{figure}[htbp]
\centering
\includegraphics[width=\textwidth]{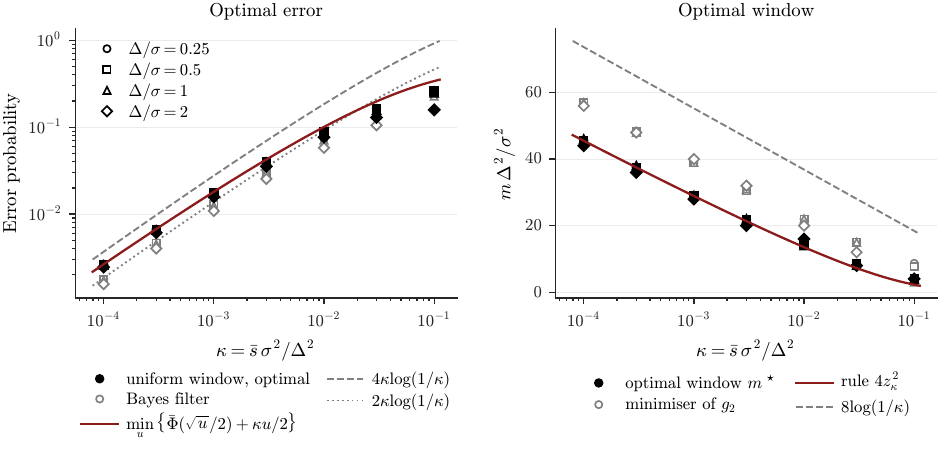}
\caption{Experiment~I, Gaussian residuals. Optimal error (left) and optimal window in resolution times (right) against \(\kappa\), with one marker shape for each \(\Delta/\sigma\).}
\label{fig:expII-collapse}
\end{figure}

Measured in resolution times, the optimal window and its error are essentially functions of \(\kappa\) alone (Figure~\ref{fig:expII-collapse}). At \(\kappa=\G{gauss-collapse-k}\) the optimal error is between \(\G{gauss-collapse-lo}\) and \(\G{gauss-collapse-hi}\) for \(\Delta/\sigma\le1\) and \(\G{gauss-collapse-d2}\) for \(\Delta/\sigma=2\), where windows are only a few observations long, and the optimal window is \(\G{gauss-collapse-m-lo}\)--\(\G{gauss-collapse-m-hi}\) resolution times for all four. Both bounds of the half-window law hold in every configuration, including asymmetric switching and unequal noise levels, and at the optimum the factor \(p_{\min}^{m-1}\) between them is above \(1-\G{gauss-sandwich-gap-max}\%\) whenever \(\kappa\le10^{-3}\). At equal \(\kappa\le0.03\), asymmetric switching changes the optimal error by at most \(\G{gauss-skew-pstar-dev}\%\) and the optimal window by at most \(\G{gauss-skew-mstar-dev}\%\). The minimizer \(m_A\) of the law is within \(\G{gauss-rule-excess-max}\%\) of the optimal error for \(\kappa\le0.03\), against \(\G{gauss-g2-excess-lo}\)--\(\G{gauss-g2-excess-hi}\%\) for the minimizer of \(g_2\). The asymptotics of Corollary~\ref{cor:window-rule}(iii) are approached slowly, as usual for \(\log(1/\kappa)\) rates. At \(\kappa=10^{-4}\) the optimal error is \(\G{gauss-asym-ratio-lo}\)--\(\G{gauss-asym-ratio-hi}\) times \(4\kappa\log(1/\kappa)\), and the filter error is \(\G{gauss-filt-ratio-k4}\) times \(2\kappa\log(1/\kappa)\).

\begin{figure}[htbp]
\centering
\includegraphics[width=\textwidth]{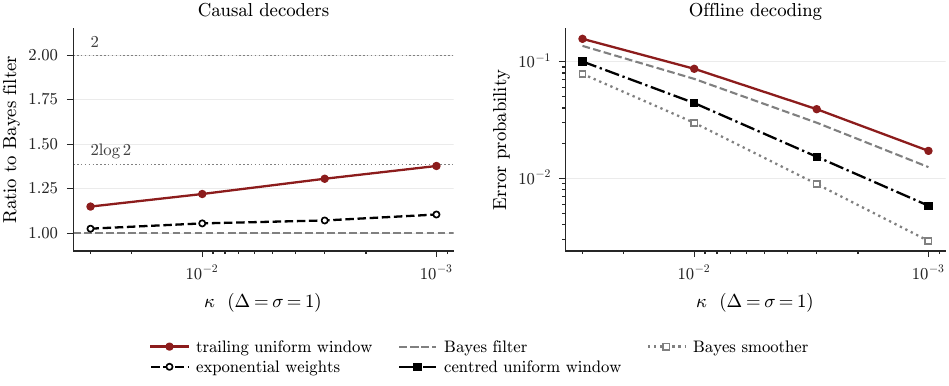}
\caption{Experiment~I, \(\Delta=\sigma=1\). Left: best uniform window and best exponential weights relative to the Bayes filter. Right: offline decoding.}
\label{fig:expII-decoders}
\end{figure}

Figure~\ref{fig:expII-decoders} compares window shapes. For \(\kappa\) between \(10^{-3}\) and \(0.03\), exponential weights are within \(\G{gauss-sim-ratio-e-lo}\)--\(\G{gauss-sim-ratio-e-hi}\) times the error of the Bayes filter (standard errors at most \(\G{gauss-sim-ratio-e-se-max}\)) and uniform windows within \(\G{gauss-sim-ratio-u-lo}\)--\(\G{gauss-sim-ratio-u-hi}\) times. Both ratios grow towards the limits \(2\log2\) and \(2\) of Proposition~\ref{prop:shape} as \(\kappa\) decreases. Offline, the centred window reduces the error of the trailing window by a factor between \(\G{gauss-offline-lo}\) and \(\G{gauss-offline-hi}\), and the Bayes smoother by a further factor between \(\G{gauss-smoother-gain-lo}\) and \(\G{gauss-smoother-gain-hi}\).

\emph{Non-Markov regimes.} We test Theorem~\ref{thm:half-window}(i) with Gaussian residuals, \(\Delta=\sigma=1\) and sojourns of mean \(\G{gauss-nm-mean}\) drawn from three laws, namely geometric sojourns (the Markov case), regular sojourns, each a sum of four geometric variables, and bursty sojourns, each drawn from an equal mixture of geometric variables with means \(10\) and \(190\). On \(\G{gauss-nm-R}\) trajectories of length \(\G{gauss-nm-T}\) we compare the uniform window, with the rule \(m_A=\G{gauss-nm-mA}\) that depends only on the switch rate, the hidden Markov filter with the matching switch rate, and exponential weights tuned on separate trajectories with Markov regimes.

\begin{figure}[htbp]
\centering
\includegraphics[width=\textwidth]{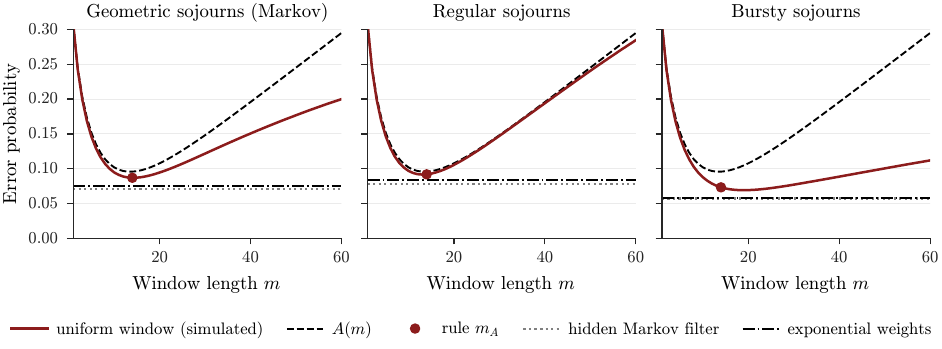}
\caption{Experiment~I, \(\Delta=\sigma=1\), with sojourns of mean \(\G{gauss-nm-mean}\) drawn from three laws.}
\label{fig:expII-nonmarkov}
\end{figure}

The bound holds for all three laws, up to simulation error (Figure~\ref{fig:expII-nonmarkov}, standard errors at most \(\G{gauss-nm-se-max}\)), and is nearly attained for regular sojourns, for which a window rarely contains two switches. At \(m_A\) the uniform window errs \(\G{gauss-nm-geo}\), \(\G{gauss-nm-reg}\) and \(\G{gauss-nm-bur}\) for geometric, regular and bursty sojourns, below \(A(m_A)=\G{gauss-nm-A}\) in each case. For bursty sojourns the best window, \(\G{gauss-nm-best-bur}\), is longer than \(m_A\). The filter remains better by a factor between \(\G{gauss-nm-ratio-lo}\) and \(\G{gauss-nm-ratio-hi}\) (standard errors at most \(\G{gauss-nm-ratio-se-max}\)), but only the window keeps a guarantee without a model of the sojourns.

Appendix~\ref{app:experiments} shows the error as a function of the age of the current regime, which isolates the half-window effect, and repeats the analysis for \(\pm1\) observations.

\subsection{Experiment II: more than two regimes}\label{sec:exp-k}

We compute the exact error of the uniform window for three regimes with means \(-1\), \(0\) and \(1\) and \(\sigma=1\), so that \(\Delta=1\), by a backward recursion over the current regime and the composition of the window. Three chains have the stationary law \((\tfrac13,\tfrac13,\tfrac13)\) and the same switch rate, namely a birth--death chain, which switches only between neighbours, a chain that switches to each of the other two regimes with equal probability, so that a third of the switches jump over the middle regime, and the cycle \(1\to2\to3\to1\), which is not reversible and has the same proportion of jumps. The mean cost of a switch is \(\bar c=1/2\) for the first chain and \(7/12\) for the other two. We take \(\kappa\) between \(\G{kreg-kappa-lo}\) and \(\G{kreg-kappa-hi}\), and simulate the Bayes filter on \(\G{kreg-Rsim}\) trajectories of length \(\G{kreg-Tsim}\).

\begin{figure}[htbp]
\centering
\includegraphics[width=\textwidth]{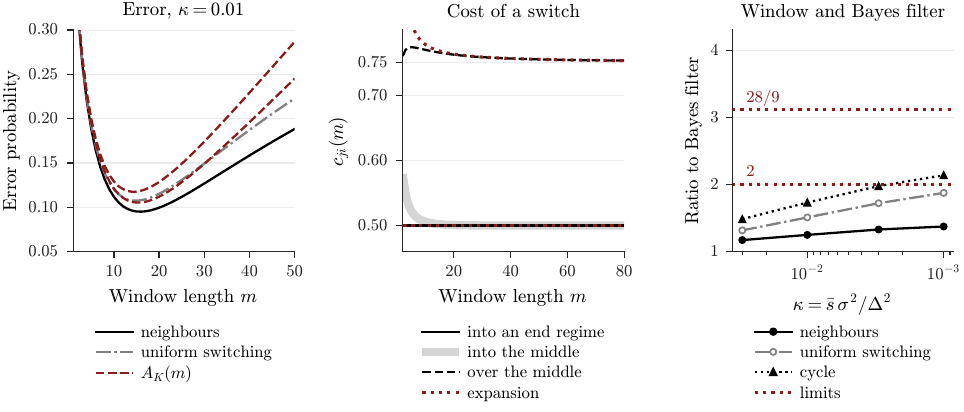}
\caption{Experiment~II, means \(-1\), \(0\), \(1\) and \(\sigma=1\). Left: error at \(\kappa=\G{kreg-kappa}\). Middle: cost of a switch (Proposition~\ref{prop:switch-geometry}). Right: ratio of the optimal window error to the Bayes filter error, with the limits of Corollary~\ref{cor:k-window}.}
\label{fig:expII-kregimes}
\end{figure}

Both bounds of the switch-cost law hold in all \(\G{kreg-nconfig}\) configurations, and the costs are close to their limits at practical windows (Figure~\ref{fig:expII-kregimes}). At the optimal window they lie between \(0.500\) and \(\G{kreg-cost-adj-max}\) for switches between neighbours and between \(\G{kreg-cost-jump-lo}\) and \(\G{kreg-cost-jump-hi}\) for jumps, and the expansion of Proposition~\ref{prop:switch-geometry} is off by at most \(\G{kreg-em-40}\) at \(m=40\).

At \(\kappa=\G{kreg-kappa}\) the optimal error is \(\G{kreg-bd-pstar}\) with switches between neighbours and \(\G{kreg-uni-pstar}\) with uniform switching, \(\G{kreg-jump-excess}\%\) more at the same switch rate, whereas the two-regime law \(A(m)\), which depends on the chain only through \(\bar s\), is the same for both. Jumps raise the error, but the minimizers of \(A_K\) and of \(A\) are within \(\G{kreg-rule-A-excess}\%\) and \(\G{kreg-two-law-excess-hi}\%\) of the optimal error in every configuration. At \(\kappa=10^{-4}\) the optimal error is \(\G{kreg-asym-lo}\)--\(\G{kreg-asym-hi}\) times \(8\bar c\,\kappa\log(1/\kappa)\), as slow an approach as for two regimes.

At \(\kappa=10^{-3}\) the window errs \(\G{kreg-win-bd-k3}\) with switches between neighbours and \(\G{kreg-win-uni-k3}\) with uniform switching, the filter \(\G{kreg-filt-bd-k3}\) and \(\G{kreg-filt-uni-k3}\), so the ratio of the window error to the filter error rises from \(\G{kreg-winratio-bd-k3}\) to \(\G{kreg-winratio-uni-k3}\), and to \(\G{kreg-winratio-cyc-k3}\) for the cycle (standard errors at most \(\G{kreg-winratio-se-max}\)), whose filter also knows which regime comes next. The simulated filter errors are \(\G{kreg-kz-lo}\)--\(\G{kreg-kz-hi}\) times the asymptotic formula \eqref{eq:kz-K} (standard errors at most \(\G{kreg-kz-se-max}\)), including the birth--death chain and the cycle, which forbid some transitions. With unequally spaced means \(-1\), \(0\) and \(2\), reversing the direction of the cycle changes the optimal error at \(\kappa=0.01\) from \(\G{kreg-rev-fwd}\) to \(\G{kreg-rev-rev}\), a factor of \(\G{kreg-rev-ratio}\), against the asymptotic factor \(\G{kreg-rev-limit}\) given by the mean costs.

\subsection{Experiment III: the reflected walk}\label{sec:exp-refl}

We simulate model~(c) with \((p,\Delta)=(\G{refl-p},\G{refl-delta})\) on lattices of size \(M\) between \(\G{refl-M-lo}\) and \(\G{refl-M-hi}\), with holding and with bouncing walls, from the stationary law of \((X_t,E_t)\), on \(\G{refl-R}\) trajectories of length \(\G{refl-T}\) (standard errors at most \(\G{refl-se-max}\)), and compare the decoder on raw increments with the position-aware decoder \eqref{eq:position-aware}.

\begin{figure}[htbp]
\centering
\includegraphics[width=\textwidth]{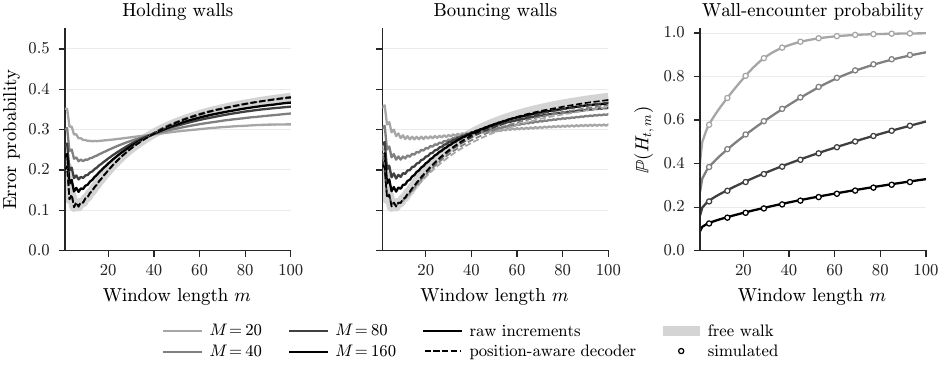}
\caption{Experiment~III, \((p,\Delta)=(\G{refl-p},\G{refl-delta})\). Error with holding (left) and bouncing (middle) walls, and the probability that a window meets a wall (right).}
\label{fig:expIII-reflected}
\end{figure}

With raw increments and holding walls (Figure~\ref{fig:expIII-reflected}), at \(m=\G{refl-m5}\) the error is \(\G{refl-err5-lo}\) for \(M=\G{refl-M-hi}\) and \(\G{refl-err5-hi}\) for \(M=\G{refl-M-lo}\), against \(\G{refl-err5-proxy}\) for the free walk, and the optimal windows, \G{refl-mstar-list} for increasing \(M\), are longer than the \(\G{refl-proxy-mstar}\) of the free walk. The walk spends between \(\G{refl-wall-lo}\%\) and \(\G{refl-wall-hi}\%\) of the time at a wall, and a window of length \(\G{refl-m5}\) meets one with probability between \(\G{refl-ph5-lo}\) and \(\G{refl-ph5-hi}\). The error of the position-aware decoder equals that of the free walk to within \(\G{refl-pa-maxdev}\) at every lattice size and window length, as Section~\ref{sec:reflected} predicts.

With bouncing walls, between \(\G{refl-b-wall-lo}\%\) (\(M=\G{refl-M-hi}\)) and \(\G{refl-b-wall-hi}\%\) (\(M=\G{refl-M-lo}\)) of the steps start at a wall and carry no information. The raw-increment decoder loses accuracy as with holding walls, with a best error of \(\G{refl-b-raw-hi}\) for \(M=\G{refl-M-lo}\). The best error of the position-aware decoder is \(\G{refl-b-best-hi}\) for \(M=\G{refl-M-lo}\) and \(\G{refl-b-best-lo}\) for \(M=\G{refl-M-hi}\), against \(\G{refl-proxy-best}\). Thinning every window by the fraction of wall steps would predict \(\G{refl-b-pred-hi}\) for \(M=\G{refl-M-lo}\). The smaller loss is not explained by Proposition~\ref{prop:position-aware}, and we leave it open.

\FloatBarrier
\section{Discussion}\label{sec:conclusion}

The results quantify the loss incurred by using a window decoder instead of the Bayes filter. A window uses only the regime means, and its length follows from the noise level and the switch rate. For two regimes the half-window bound holds for every stationary switching law, so it remains applicable when the sojourn distribution is unknown, as for bull and bear markets \cite{LundeTimmermann2004}. Under rare Markov switching and Gaussian noise with common variance, the comparison with the Bayes filter, which needs the full model, is sharp, and the uniform window loses asymptotically a factor two.

The same comparison applies to agents that act on the decoded regime, since their avoidable regret is proportional to the decoding error (Remark~\ref{rem:regret}). Linear recurrent memories perform well in partially observable reinforcement learning \cite{ZhaoEtAl2026}. For two symmetric regimes with rare switches and Gaussian noise with common variance the results measure what such a memory gives up against the belief state, asymptotically at most a factor \(2\log2\) for exponential weights, and Proposition~\ref{prop:two-block} shows that a linear memory whose weights depend on the switch rate attains the Bayes error to first order.

A single window cannot adapt to small and large switches at once. The simple filters of Khasminskii and Zeitouni \cite{KhasminskiiZeitouni1996} close this gap asymptotically with a coarse test followed by likelihood tests over windows whose lengths depend on the pair of regimes. For Gaussian residuals with a common variance these tests reduce to window sums. This suggests that a decoder that combines windows of several lengths, with no more information than the window rule, may attain the Bayes rate, and a non-asymptotic analysis of such a decoder is open. For two regimes a single linear decoder already attains it to first order (Proposition~\ref{prop:two-block}), with a bound that exceeds the leading term of the Bayes error by a relative amount of order \((\log(1/\kappa))^{-1/4}\). We do not know whether a faster rate is possible.

For two regimes, the lower bound in Theorem~\ref{thm:half-window} relies on the Markov property, and Experiment~I shows that this matters. Under bursty sojourns with the same switch rate, windows frequently contain several switches, the error is substantially smaller than \(A(m)\), and the optimal window is longer than \(m_A\). Which features of the switching law, beyond the switch rate, determine the error? A sharp answer would require the joint distribution of the switches within a window.

The window rule requires the regime means, the noise level and the switch rate. How long a single trajectory must be to estimate them, including from labels obtained in hindsight, is a separate question. Without labels, the means must be estimated by clustering, the expectation--maximization algorithm or spectral methods, and the learnability thresholds of \cite{AbrahamGassiatNaulet2023} become relevant.

When the regimes differ in diffusivity rather than in drift, as in single-particle tracking \cite{PerssonEtAl2013}, the decoder must be based on a feature such as the squared increments. For centred Gaussian increments with variances \(\sigma_i^2\), the squared increments are scaled \(\chi^2_1\) variables, which are stochastically ordered and have exponential moments near zero, so Theorem~\ref{thm:half-window} and Corollary~\ref{cor:window-rule}(iii) apply to a uniform window of squared increments. For other increment laws these two properties must be checked. The margin bounds would require sub-exponential residuals, and how best to combine drift and diffusivity features is open.

For the reflected walk with bouncing walls, the position-aware decoder loses much less than the fraction of uninformative wall steps suggests (Experiment~III), and we cannot yet explain why. The distribution of the information clock would control one term of Proposition~\ref{prop:position-aware}, but it would not by itself explain the classification error. Finally, when the actions determine what is observed, as in regime-switching bandits \cite{ZhouXiongChenGao2021}, exploration interacts with recovery, and the regret identity of Remark~\ref{rem:regret} no longer applies directly.


\clearpage
\appendix

\section{Proofs}\label{app:proofs}

\subsection{Proofs for Section~\ref{sec:half-window}}\label{app:proofs-local}

\begin{proof}[Proof of Proposition~\ref{thm:recoverability}]
Let \(S_{t,m}:=\sum_{s=t-m+1}^t\xi_s\). On \(\mathrm{NoSwitch}_{t,m}\), \(\bar Y_{t,m}-\mu_{E_t}=S_{t,m}/m\), and if \(|S_{t,m}|<m\Delta/2\), then \(\mu_{E_t}\) is the unique nearest mean, since for \(j\neq E_t\)
\[
|\bar Y_{t,m}-\mu_j|\ge|\mu_j-\mu_{E_t}|-|S_{t,m}|/m>\Delta/2 .
\]
Hence
\[
\PP(\widehat E_t\neq E_t)\le\PP\bigl(|S_{t,m}|\ge m\Delta/2\bigr)+\PP(\mathrm{NoSwitch}_{t,m}^c).
\]
For \(s<t\), \(\xi_s\) is \(\mathcal H_t\)-measurable, so conditioning successively on \(\mathcal H_t,\mathcal H_{t-1},\ldots\) and using Assumption~\ref{ass:residual-model}(ii) gives \(\EE[e^{\lambda S_{t,m}}]\le e^{m\sigma^2\lambda^2/2}\) for every \(\lambda\in\RR\), the argument behind the Azuma--Hoeffding inequality \cite{Azuma1967}. The Chernoff bound on both tails gives \(\PP(|S_{t,m}|\ge m\Delta/2)\le2e^{-m\Delta^2/(8\sigma^2)}\), and \eqref{eq:noswitch-exact} completes the proof.
\end{proof}

\emph{The margin of a window.} The proofs of Propositions~\ref{prop:shape} and~\ref{prop:centred} use a bound for general linear decoders under martingale-difference residuals. Let \(K=2\), and consider the linear decoder \eqref{eq:linear-decoder} with weights \(w\) on a finite index set \(\mathcal T\), which may contain times after \(t\). The weight carried by the current regime and the resulting margin are
\begin{equation}
W_t:=\sum_{s\in\mathcal T}w_s\mathbf 1\{E_s=E_t\},
\qquad
\frac{\Delta}{2}\,(2W_t-1).
\label{eq:margin}
\end{equation}
The margin is the distance from \(\theta\) to the conditional mean of \(S_t\) given the regimes, signed so that it is positive when the decoder is correct on average.

\begin{lemma}[Moment generating function along a regime path]
\label{lem:path-mgf}
Let Assumption~\ref{ass:residual-model} hold, let \(1\le a\le b\), let \(e=(e_a,\ldots,e_b)\) be a regime path, and let \(c_s\in\RR\) for \(s\in\mathcal T\subset\{a,\ldots,b\}\). Then
\[
\EE\Bigl[\mathbf 1\{E_{a:b}=e\}\exp\Bigl(\sum_{s\in\mathcal T}c_s\xi_s\Bigr)\Bigr]
\le
\exp\Bigl(\frac{\sigma^2}{2}\sum_{s\in\mathcal T}c_s^2\Bigr)\,\PP(E_{a:b}=e).
\]
\end{lemma}

\begin{proof}
Put \(c_s:=0\) for \(s\notin\mathcal T\) and \(Z_r:=\mathbf 1\{E_{a:r}=e_{a:r}\}\exp(\sum_{s=a}^r c_s\xi_s)\) for \(a\le r\le b\). The variable \(Z_{r-1}\) is measurable with respect to \(\mathcal F_{r-1}:=\sigma(E_1,\ldots,E_{r-1},Y_1,\ldots,Y_{r-1})\subset\mathcal H_r\), and \(E_r\) is \(\mathcal H_r\)-measurable. Assumption~\ref{ass:residual-model}(ii) and then~(i) give
\[
\EE[Z_r]
\le e^{\sigma^2c_r^2/2}\,\EE\bigl[Z_{r-1}\mathbf 1\{E_r=e_r\}\bigr]
= e^{\sigma^2c_r^2/2}\,Q_{e_{r-1}e_r}\,\EE[Z_{r-1}],
\]
because \(\EE[\mathbf 1\{E_r=e_r\}\mid\mathcal F_{r-1}]=Q_{E_{r-1}e_r}\) and \(E_{r-1}=e_{r-1}\) wherever \(Z_{r-1}\neq0\). Since \(\EE[Z_a]\le e^{\sigma^2c_a^2/2}\PP(E_a=e_a)\), induction gives the claim.
\end{proof}

\begin{theorem}[Margin bound]
\label{thm:margin}
Under Assumption~\ref{ass:residual-model}, the linear decoder \eqref{eq:linear-decoder} satisfies, counting ties as errors,
\begin{equation}
\PP(\widehat E_t\neq E_t)
\le
\EE\!\left[\exp\!\left(-\frac{\Delta^2\,(2W_t-1)_+^2}{8\sigma^2\|w\|^2}\right)\right],
\qquad
\|w\|^2:=\sum_{s\in\mathcal T}w_s^2 .
\label{eq:margin-bound}
\end{equation}
\end{theorem}

\begin{proof}
Let \([a,b]\) be the smallest interval containing \(\mathcal T\cup\{t\}\) and fix a regime path \(e\) on it, with margin \(M(e)=\frac{\Delta}{2}(2W(e)-1)\). Write \(\epsilon:=1\) if \(e_t=1\) and \(\epsilon:=-1\) if \(e_t=2\). Since \(\mu_{e_s}-\theta=\pm\Delta/2\) according to whether \(e_s\) is regime~2 or~1, on \(\{E_{a:b}=e\}\) the decoder errs if and only if \(\epsilon\sum_{s\in\mathcal T}w_s\xi_s\ge M(e)\). For \(\lambda\ge0\), Markov's inequality and Lemma~\ref{lem:path-mgf} give
\[
\PP(\widehat E_t\neq E_t,\ E_{a:b}=e)
\le e^{-\lambda M(e)+\lambda^2\sigma^2\|w\|^2/2}\,\PP(E_{a:b}=e).
\]
Take \(\lambda=M(e)_+/(\sigma^2\|w\|^2)\) and sum over \(e\).
\end{proof}

For the trailing uniform window, \(W_t\ge\min\{\tau_t,m\}/m\), so the bound depends only on how long ago the last switch occurred. For a stationary two-state chain the age has a simple law. A switch into regime~\(i\) at a given time has probability \(\pi_jQ_{ji}=\pi_iQ_{ij}=\bar s/2\), \(j\neq i\), so that
\begin{equation}
\PP(\tau_t=k,\ E_t=i)=\frac{\bar s}{2}\,Q_{ii}^{\,k-1},\qquad k\ge1 .
\label{eq:age-law}
\end{equation}

\begin{corollary}[Half-window bound]
\label{cor:half-window}
Under Assumption~\ref{ass:residual-model}, for a stationary two-state chain and the trailing uniform window of length \(m\le t\),
\begin{align}
\PP(\widehat E_t\neq E_t)
&\le
\PP(\tau_t\ge m)\,e^{-m\Delta^2/(8\sigma^2)}
+\frac{\bar s}{2}\sum_{k=1}^{m-1}\bigl(Q_{11}^{k-1}+Q_{22}^{k-1}\bigr)\exp\!\left(-\frac{\Delta^2\,(2k-m)_+^2}{8\sigma^2m}\right)
\notag\\
&\le
e^{-m\Delta^2/(8\sigma^2)}+\bar s\left(\frac m2+1+\frac{\sigma}{\Delta}\sqrt{\frac{\pi m}{2}}\right).
\label{eq:half-window-bound}
\end{align}
\end{corollary}

\begin{proof}
Insert \(W_t\ge\min\{\tau_t,m\}/m\) and \(\|w\|^2=1/m\) into \eqref{eq:margin-bound} and use \eqref{eq:age-law}. For the second line, bound \(Q_{ii}^{k-1}\) by one. The \(\lfloor m/2\rfloor\) terms with \(2k\le m\) are at most one, and the remaining terms, indexed by \(j=2k-m\ge1\) in steps of two, sum to at most
\[
1+\frac12\int_0^\infty e^{-x^2\Delta^2/(8\sigma^2m)}\,dx=1+\frac{\sigma}{\Delta}\sqrt{\frac{\pi m}{2}} .\qedhere
\]
\end{proof}

The switching term of \(g_2\) is \(1-p_{\min}^{m-1}\approx(1-p_{\min})\,m\), whereas here it is about \(\bar s\,m/2\). Under martingale-difference residuals, a switch thus causes an error essentially only when it lies in the most recent half of the window, which Figure~\ref{fig:switch-age} shows directly. Under independence, Theorem~\ref{thm:half-window} makes the average cost of a switch exactly one half.

\emph{Unequal noise levels.} Gaussian residuals with standard deviations \(\sigma_1\neq\sigma_2\) are not stochastically ordered, but \(\psi_1\) and \(\psi_2\) are still nonincreasing, which is all the proof of Theorem~\ref{thm:half-window} uses. With \(k\) observations from the current regime~1, \(\psi_1(k)=\Phibar(f(k))\) with \(f(k):=\Delta(k-m/2)/(k\sigma_1^2+(m-k)\sigma_2^2)^{1/2}\), whose derivative has the sign of \(k(\sigma_1^2-\sigma_2^2)/2+m(\sigma_1^2+3\sigma_2^2)/4\). This is affine in \(k\) and positive at \(k=0\) and \(k=m\), hence on \([0,m]\), and symmetrically for \(\psi_2\).

\begin{corollary}[Sharper bounds for the symmetric chain]
\label{cor:symmetric-bounds}
In the setting of Theorem~\ref{thm:half-window}(i), with \(F_1\) stochastically smaller than \(F_2\), let \((E_t)\) be the symmetric chain \(Q(p)\) and \(r:=\lceil m/2\rceil-1\). Then
\begin{equation}
\bar P(m)+(m-1)(1-p)p^{m-2}\bigl(\tfrac12-\bar P(m)\bigr)
\le\PP(\widehat E_t\neq E_t)\le
p^{m-1}\bar P(m)+1-p^{r}+\mathbf 1\{m\text{ even}\}\,\tfrac12(1-p)p^{m/2-1}.
\label{eq:symmetric-bounds}
\end{equation}
\end{corollary}

The upper bound is at most the single-regime error plus the probability \(1-p^{r}\) that the last switch lies in the more recent half of the window, with half weight for a switch exactly in the middle, and it never exceeds one. The lower bound is at least the single-regime error whenever \(\bar P(m)\le\tfrac12\).

\begin{proof}[Proof of Corollary~\ref{cor:symmetric-bounds}]
For the symmetric chain, \(\PP(\tau_t=k,E_t=i)=\tfrac12(1-p)p^{k-1}\) for \(k\ge1\), and \(\PP(\tau_t\ge m,E_t=i)=\tfrac12p^{m-1}\). The argument of Theorem~\ref{thm:half-window}(i) therefore gives
\[
\PP(\widehat E_t\neq E_t)\le p^{m-1}\bar P(m)+\frac{1-p}{2}\sum_{k=1}^{m-1}p^{k-1}\bigl(\psi_1(k)+\psi_2(k)\bigr).
\]
For \(k<m/2\), pairing gives \(\psi_2(m-k)=1-\psi_1(k)\), so
\[
p^{k-1}\psi_1(k)+p^{m-k-1}\psi_2(m-k)\le p^{k-1}\bigl(\psi_1(k)+\psi_2(m-k)\bigr)=p^{k-1},
\]
and likewise with the regimes exchanged. For even \(m\) the term \(k=m/2\) contributes \(p^{m/2-1}\bigl(\psi_1(m/2)+\psi_2(m/2)\bigr)=p^{m/2-1}\). Summing over the \(r\) values \(k<m/2\),
\[
\frac{1-p}{2}\sum_{k=1}^{m-1}p^{k-1}\bigl(\psi_1(k)+\psi_2(k)\bigr)
\le(1-p)\sum_{k=1}^{r}p^{k-1}+\mathbf 1\{m\text{ even}\}\,\frac{1-p}{2}\,p^{m/2-1},
\]
which is the upper bound. For the lower bound, keep the events \(\mathcal A_{k,i}\) of the proof of Theorem~\ref{thm:half-window}(ii), which have probability \(\tfrac12(1-p)p^{m-2}\) and error \(\psi_i(k)\). On the rest of \(\{\tau_t=k,\ E_t=i\}\), which has probability \(\tfrac12(1-p)(p^{k-1}-p^{m-2})\), the error is at least \(\psi_i(m)=P^{(i)}(m)\) by monotonicity. Since \(\sum_{k<m}(1-p)p^{k-1}=1-p^{m-1}\),
\[
\PP(\widehat E_t\neq E_t)\ge\bar P(m)-(m-1)(1-p)p^{m-2}\bar P(m)+\frac{(1-p)p^{m-2}}{2}\sum_{k=1}^{m-1}\bigl(\psi_1(k)+\psi_2(k)\bigr),
\]
and the pairing identity gives the lower bound.
\end{proof}

\emph{Comparing all windows.} The sharp asymptotics of Corollary~\ref{cor:window-rule}, Proposition~\ref{prop:shape} and Corollary~\ref{cor:k-window} rest on one lemma, which bounds the error of every competing window from below at once. What makes the comparison uniform is that the lower bound~(iii) holds for every window length, with functions \(F\) and \(r\) that do not depend on the switching probabilities.

\begin{lemma}[Comparing all windows]
\label{lem:three-ranges}
Let \(s\to0\), write \(\ell:=\log(1/s)\), and for each \(s\) let \(P_s\ge0\) be a function on a set \(\mathcal M_s\subset[1,C/s]\) of window lengths, with \(C\) fixed. Let \(a,b,I>0\) be constants and \(F,r\) functions on \([1,\infty)\), not depending on \(s\), such that
\begin{enumerate}[label=(\roman*)]
\item for every \(\delta>0\) there is \(c_\delta>0\) with \(F(m)\ge c_\delta e^{-m(I+\delta)}\) for all \(m\ge1\),
\item \(r\ge0\), and \(r(m)\to1\) as \(m\to\infty\),
\item for all \(m\in\mathcal M_s\),
\[
P_s(m)\ge(1-as)^m\bigl[F(m)+b\,s\,(m-1)\,r(m)\bigr]-\epsilon_s,
\]
where \(\epsilon_s=o(s\ell)\) does not depend on \(m\).
\end{enumerate}
Then \(\inf_{m\in\mathcal M_s}P_s(m)\ge(bs\ell/I)(1-o(1))\). If moreover \(m_s\in\mathcal M_s\) and \(P_s(m_s)\le(bs\ell/I)(1+o(1))\), then \(m_s=(\ell/I)(1+o(1))\).
\end{lemma}

\begin{proof}
Fix \(\eta\in(0,1/2)\), put \(\delta:=\eta I/2\), and let \(s\) be small. We use \((1-as)^m\ge1-asm\) and, for \(as\le1/2\), \((1-as)^m\ge e^{-2asm}\). Split \(\mathcal M_s\) into three ranges.
\begin{enumerate}[label=(\alph*)]
\item \(m\le(1-\eta)\ell/I\). Then \((1-as)^m\ge1-as\ell/I\ge1/2\), and since \((1-\eta)(1+\eta/2)\le1-\eta/2\),
\[
P_s(m)\ge\tfrac12c_\delta\,s^{(1-\eta)(1+\eta/2)}-\epsilon_s\ge\tfrac12c_\delta\,s^{1-\eta/2}-\epsilon_s .
\]
\item \((1-\eta)\ell/I<m\le1/(as\ell)\). Then \((1-as)^m\ge1-1/\ell\). Since \(r\) does not depend on \(s\) and \(r(m)\to1\), \(r\ge1-\eta\) on \([(1-\eta)\ell/I,\infty)\) for small \(s\), and
\[
P_s(m)\ge\Bigl(1-\frac1\ell\Bigr)(1-\eta)\,b\,s\,(m-1)-\epsilon_s\ge(1-\eta)^2\,\frac{bs\ell}{I}\,(1-o(1)).
\]
\item \(1/(as\ell)<m\le C/s\). Then \((1-as)^m\ge e^{-2aC}\) and, for small \(s\), \(r(m)\ge1/2\), so
\[
P_s(m)\ge e^{-2aC}\,\frac b2\,\Bigl(\frac1{a\ell}-s\Bigr)-\epsilon_s .
\]
\end{enumerate}
Each bound holds uniformly over its range, and the bounds in (a) and (c), divided by \(s\ell\), tend to infinity. Hence \(\inf_{\mathcal M_s}P_s\ge(1-\eta)^2(bs\ell/I)(1-o(1))\), and \(\eta\) is arbitrary. For the second claim, (a) and (c) show that \(m_s\) lies in range~(b) for small \(s\), so \(m_s>(1-\eta)\ell/I\), and by~(iii)
\[
\Bigl(1-\frac1\ell\Bigr)r(m_s)\,b\,s\,(m_s-1)-\epsilon_s\le P_s(m_s)\le\frac{bs\ell}{I}\,(1+o(1)).
\]
Since \(m_s\to\infty\), \(r(m_s)\to1\), so \(m_s\le(\ell/I)(1+o(1))\), and \(\eta\) is arbitrary.
\end{proof}

\emph{Long windows.} Lemma~\ref{lem:three-ranges} compares windows of length at most \(C/s\). For two regimes, longer windows are excluded by the identity~\eqref{eq:reversal}. For more regimes the chain need not be reversible, and a rotation of sojourns replaces the time reversal.

\begin{lemma}[Long windows]
\label{lem:long-windows}
Let \((E_t)\) be a stationary Markov chain with transition matrix \(Q\), let the observations be independent given the regime path with \(Y_s\sim F_{E_s}\), and apply the decoder \eqref{eq:decoder} to the trailing window of length \(m\le t\). Let
\[
\rho:=\min\Bigl\{1,\ \frac{\pi_bQ_{ca}}{\pi_aQ_{ab}}\ :\ a\neq b,\ c\neq a,\ Q_{ab}>0,\ Q_{ca}>0\Bigr\}.
\]
Then
\[
\PP(\widehat E_t\neq E_t)\ge\frac{\rho}{1+\rho}\,\PP\bigl(E_{t-m+1}\neq E_t,\ Q_{E_tE_{t-m+1}}>0\bigr).
\]
For two regimes \(\rho=1\), and the bound reduces to the switching term \(\tfrac12\PP(E_{t-m+1}\neq E_t)\) of Theorem~\ref{thm:half-window}(ii).
\end{lemma}

\begin{proof}
Given the regime path \(e\) of the window, the window sum has a law that depends only on the composition \(n(e)\) of the window, the number of its observations from each regime, so the decoder returns regime~\(r\) with a probability \(g_r(n(e))\). Writing \(\PP(n)\) and \(\PP(n,r)\) for the probabilities that the window has composition \(n\), and in addition that \(E_t=r\),
\[
\PP(\widehat E_t\neq E_t)=\sum_n\Bigl[\PP(n)-\sum_rg_r(n)\,\PP(n,r)\Bigr]\ge\sum_n\Bigl[\PP(n)-\max_r\PP(n,r)\Bigr].
\]
Write a path as a sequence of sojourns in regimes \(r_1,\ldots,r_J\) with lengths \(\ell_1,\ldots,\ell_J\), oldest first, so that
\[
\PP(e)=\pi_{r_1}\prod_{j=1}^JQ_{r_jr_j}^{\ell_j-1}\prod_{j=1}^{J-1}Q_{r_jr_{j+1}},
\]
and let \(B\) be the set of paths with \(r_J\neq r_1\) and \(Q_{r_Jr_1}>0\). Moving the oldest sojourn to the end maps \(e\in B\) to the path \(e'\) with sojourns \((r_2,\ell_2),\ldots,(r_J,\ell_J),(r_1,\ell_1)\). The map is injective, \(e'\) has the same composition as \(e\) and current regime \(r_1\neq r_J\), and
\[
\frac{\PP(e')}{\PP(e)}=\frac{\pi_{r_2}Q_{r_Jr_1}}{\pi_{r_1}Q_{r_1r_2}}\ge\rho .
\]
Fix \(n\) and let \(r^\star\) maximize \(\PP(n,r)\). The map sends the paths of \(B\) with composition \(n\) and current regime \(r^\star\) to paths with composition \(n\) and another current regime, so \(\PP(n)-\PP(n,r^\star)\ge\rho\,\PP(B,n,r^\star)\), and trivially \(\PP(n)-\PP(n,r^\star)\ge\PP(B,n)-\PP(B,n,r^\star)\). Multiplying the second inequality by \(\rho\le1\) and adding gives \(\PP(n)-\PP(n,r^\star)\ge\rho\,\PP(B,n)/(1+\rho)\), and summing over \(n\) proves the bound. For two regimes \(r_2=r_J\) on \(B\), so the ratio is \(\pi_{r_2}Q_{r_2r_1}/(\pi_{r_1}Q_{r_1r_2})=1\), and \(\PP(E_{t-m+1}\neq E_t)=2\pi_1\pi_2(1-\lambda^{m-1})\).
\end{proof}

\begin{proof}[Proof of Corollary~\ref{cor:window-rule}]
(i) Let \(p_{\min}>0\), and write \(P(m)\) for the error of the window of length \(m\). For every \(m\ge1\), the two bounds of Theorem~\ref{thm:half-window}(ii) and the minimality of \(m_A\) give \(P(m_A)\le A(m_A)\le A(m)\le P(m)/p_{\min}^{\,m-1}\). If the infimum \(P^\star\) is attained at \(m^\star\), taking \(m=m^\star\) gives the bound \(P^\star/p_{\min}^{\,m^\star-1}\).

(ii) Setting the derivative to zero gives \(\varphi(z)\Delta^2/(4\sigma^2z)=\bar s\) with \(z=\sqrt m\,\Delta/(2\sigma)\).

(iii) Write \(\varepsilon:=1-p_{\min}\), \(\ell:=\log(1/\bar s)\), and \(P(m)\) for the error of the window of length \(m\). Since the switching probabilities have a fixed ratio, \(\pi\) and \(a:=\varepsilon/\bar s\) are constants, and \(a\ge1\) because \(\bar s\le\varepsilon\).

\emph{Long windows.} For small \(\varepsilon\), \(0\le\lambda=1-Q_{12}-Q_{21}\le1-\varepsilon\), so Theorem~\ref{thm:half-window}(ii) gives \(P(m)\ge\pi_1\pi_2\bigl(1-(1-\varepsilon)^{m-1}\bigr)\ge\pi_1\pi_2(1-e^{-1})\) whenever \(m-1\ge1/\varepsilon\). By the upper bound below, \(P^\star\to0\), so for small \(\bar s\) the optimum lies among the windows of length at most \(2/\varepsilon\le2/\bar s\), and it suffices to apply Lemma~\ref{lem:three-ranges} with \(s:=\bar s\) to these windows.

\emph{Single-regime errors.} By Chernoff's bound, \(P^{(i)}(m)\le e^{-mI_i}\). The lower bound in Cram\'er's theorem \cite[Theorem~2.2.3]{DemboZeitouni1998}, applied to the open half-line beyond \(\theta\), and the continuity of the rate function at \(\theta\), an interior point of its domain, give \(\liminf_m m^{-1}\log P^{(i)}(m)\ge-I_i\). Since \(P^{(i)}(m)>0\) for every \(m\), for every \(\delta>0\) there is \(c_\delta>0\) with \(P^{(i)}(m)\ge c_\delta e^{-m(I_i+\delta)}\) for all \(m\). Hence, with \(\pi_{\min}:=\min_i\pi_i\),
\[
\pi_{\min}\,c_\delta\,e^{-m(I+\delta)}\le\bar P(m)\le e^{-mI}.
\]

\emph{Upper bound.} For \(m_0:=\lceil\ell/I\rceil\), which is at most \(2/\varepsilon\) for small \(\bar s\), \(\bar P(m_0)\le\bar s\), so the upper bound of Theorem~\ref{thm:half-window}(ii) gives
\[
P^\star\le P(m_0)\le A(m_0)\le\bar s+\frac{\bar s\,\ell}{2I}=\frac{\bar s\,\ell}{2I}\,(1+o(1)).
\]

\emph{Lower bound.} By Theorem~\ref{thm:half-window}(ii), for every \(m\),
\[
P(m)\ge p_{\min}^{m-1}A(m)\ge(1-a\bar s)^m\Bigl[\bar P(m)+\frac{\bar s}2\,(m-1)\Bigr],
\]
which is condition~(iii) of the lemma with \(F:=\bar P\), \(b:=1/2\), \(r\equiv1\) and \(\epsilon_s:=0\). The lemma gives \(P^\star\ge(\bar s\ell/(2I))(1-o(1))\) and, with the upper bound, \(m^\star=(\ell/I)(1+o(1))\). The function \(A\) satisfies the same condition and \(A(m_A)\le A(m_0)\). Since \(A(m)\ge\bar s(m-1)/2\), also \(m_A\le1+2A(m_0)/\bar s=O(\ell)\), and the lemma gives \(m_A=(\ell/I)(1+o(1))\).

\emph{Gaussian residuals.} Here \(I_i=(\theta-\mu_i)^2/(2\sigma^2)=\Delta^2/(8\sigma^2)\), and both \(\ell\) and \(z_\kappa^2/2\) equal \(\log(1/\kappa)(1+o(1))\).
\end{proof}

\begin{proof}[Proof of Proposition~\ref{prop:shape}]
By Cauchy--Schwarz, \(\tfrac14=(\int_0^{u_g}g)^2\le u_g\int_0^{u_g}g^2\le c_g\). Let \(G(u):=\int_0^ug\) and \(v_g:=\int g^2\), which is finite because \(g\le g(0)\). Since \(g\) is nonincreasing and \(\int_{u_g}^\infty g=1/2\), \(g\) is positive on \([0,u_g]\) and on a right neighbourhood of \(u_g\).

\emph{Weights.} By stationarity and reversibility, \(E_t,E_{t-1},\ldots,E_1\) are the first \(t\) states of a stationary symmetric chain. Theorem~\ref{thm:margin} bounds the error at time \(t\), and for Gaussian residuals the formula used below gives it, as the expectation of a bounded continuous function of these states and of the weights \(g(j/L)/\sum_{k<t}g(k/L)\), \(0\le j<t\). These converge as \(t\to\infty\), so by dominated convergence \(P_L\) obeys the same bound, and the same formula, with the weights
\[
w^{(L)}_{t-j}:=\frac{g(j/L)}{\sum_{k\ge0}g(k/L)},\qquad j\ge0,
\]
on a regime path extended to the infinite past. Write \(G_L(x):=\sum_{0\le j<x}w^{(L)}_{t-j}\). Since \(g\) is nonincreasing, comparing sums with integrals gives \(L\le\sum_{k\ge0}g(k/L)\le L+g(0)\) and \(LG(u)\le\sum_{0\le j<uL}g(j/L)\le LG(u)+g(0)\), and likewise for \(g^2\). Hence, for all \(u\ge0\) and \(L\ge1\),
\begin{equation}
|G_L(uL)-G(u)|\le\frac{g(0)}L,
\qquad
\frac{v_g}{(1+g(0)/L)^2}\le L\|w^{(L)}\|^2\le v_g+\frac{g(0)^2}L .
\label{eq:riemann}
\end{equation}

\emph{Upper bound.} Since the \(\tau_t\) most recent observations come from the current regime, \(W_t\ge G_L(\tau_t)\). Fix \(\delta\in(0,1/4)\) and let
\[
\eta:=2G(u_g+\delta)-1>0,
\qquad
u_\infty:=\inf\{u:2G(u)-1\ge1-\delta\}.
\]
By \eqref{eq:riemann}, for \(L\ge4g(0)/\min\{\eta,\delta\}\), \(2G_L(k)-1\ge\eta/2\) for \(k>(u_g+\delta)L\) and \(2G_L(k)-1\ge1-2\delta\) for \(k>u_\infty L\). Split \eqref{eq:margin-bound} according to \(\tau_t\le(u_g+\delta)L\), \((u_g+\delta)L<\tau_t\le u_\infty L\) and \(\tau_t>u_\infty L\), and use \(\PP(\tau_t=k)\le1-p\) for the symmetric chain. With \(\alpha_L:=\Delta^2/(8\sigma^2\|w^{(L)}\|^2)\),
\[
P_L\le(1-p)(u_g+\delta)L+(1-p)\,u_\infty L\,e^{-\eta^2\alpha_L/4}+e^{-(1-2\delta)^2\alpha_L}.
\]
Take \(L:=8(1+\delta)^2(1-2\delta)^{-2}v_g(\sigma^2/\Delta^2)\log(1/\kappa)\). For small \(\kappa\), \(L\|w^{(L)}\|^2\le(1+\delta)v_g\) by \eqref{eq:riemann}, so \((1-2\delta)^2\alpha_L\ge(1+\delta)\log(1/\kappa)\) and \(\alpha_L\ge\log(1/\kappa)\). The last two terms are then at most \((1-p)u_\infty L\,\kappa^{\eta^2/4}+\kappa^{1+\delta}=o(\kappa\log(1/\kappa))\), and the first is \(8(u_g+\delta)(1+\delta)^2(1-2\delta)^{-2}v_g\,\kappa\log(1/\kappa)\). Let \(\delta\to0\).

\emph{Lower bound.} Let the residuals be independent \(N(0,\sigma^2)\) variables, independent of the regimes. Given the regime path, \(S_t\) is normal with variance \(\sigma^2\|w\|^2\), and the decoder errs with probability \(\Phibar\bigl(\Delta(2W_t-1)/(2\sigma\|w\|)\bigr)\), which is at least \(\Phibar(\Delta/(2\sigma\|w\|))\) on every path. Fix \(\delta\in(0,\min\{1,u_g\})\). Since \(g>0\) on \([0,u_g]\), \(\eta:=1-2G(u_g-\delta)>0\). Choose \(c>u_g\) with \(1-G(c)\le\eta/8\). By \eqref{eq:riemann} there is \(L_0\ge(1-\delta)v_g\), depending on \(g\), \(\delta\) and \(\Delta/\sigma\) but not on \(p\), such that for every \(L\ge L_0\):
\begin{enumerate}[label=(\alph*)]
\item \(L\|w^{(L)}\|^2\ge(1-\delta)\,v_g\),
\item \(g(0)/L\le\eta/8\), so that \(G_L\bigl((u_g-\delta)L\bigr)\le G(u_g-\delta)+\eta/8\) and \(1-G_L(cL)\le\eta/4\),
\item \(\Phibar\bigl(\eta\Delta/(8\sigma\|w^{(L)}\|)\bigr)\le\delta\).
\end{enumerate}
Let \(L\ge L_0\) and, for \(1\le k\le(u_g-\delta)L\), let \(B_k\) be the event that \(\tau_t=k\) and no other switch occurs in the last \(\lceil cL\rceil\) steps. On \(B_k\) the current regime carries the weights at lags \(j<k\) and at most those at lags \(j\ge cL\), so by~(b)
\[
2W_t-1\le2\Bigl(G(u_g-\delta)+\frac\eta8+\frac\eta4\Bigr)-1=-\frac\eta4,
\]
and by~(c) the error is at least \(1-\delta\). For the symmetric chain \(\PP(B_k)\ge(1-p)p^{cL}\), and the event that no switch occurs in the last \(\lceil cL\rceil\) steps has probability at least \(p^{cL}\). These events are disjoint, so by~(a)
\begin{equation}
P_L\ge p^{cL}\Bigl[(1-\delta)(1-p)\bigl((u_g-\delta)L-1\bigr)+\Phibar\Bigl(\frac{\Delta}{2\sigma}\sqrt{\frac{L}{(1-\delta)v_g}}\Bigr)\Bigr],
\qquad L\ge L_0 .
\label{eq:shape-lower}
\end{equation}
For \(L<L_0\), \(\|w^{(L)}\|\ge w^{(L)}_t\ge g(0)/(g(0)+L_0)\), so \(P_L\ge\Phibar\bigl(\Delta(g(0)+L_0)/(2\sigma g(0))\bigr)\), a positive constant independent of \(p\).

We apply Lemma~\ref{lem:three-ranges} with \(s:=1-p\) and the variable \(m:=L/((1-\delta)v_g)\), which ranges over \(\mathcal M_s:=\{L/((1-\delta)v_g):L_0\le L\le C/(1-p)\}\subset[1,C'/s]\), \(C':=C/((1-\delta)v_g)\). Put
\[
F(m):=\Phibar\Bigl(\frac{\Delta\sqrt m}{2\sigma}\Bigr),
\qquad
b:=(1-\delta)^2(u_g-\delta)v_g,
\qquad
r\equiv1,
\qquad
a:=2c(1-\delta)v_g,
\qquad
\epsilon_s:=s .
\]
Since \(\log p\ge-2(1-p)\) for \(p\ge1/2\), \(p^{cL}\ge e^{-asm}\ge(1-as)^m\). The bracket in \eqref{eq:shape-lower} is at least \(F(m)+bs(m-1)-s\), and the Mills-ratio bound \(\Phibar(x)\ge x\varphi(x)/(1+x^2)\) gives condition~(i) of the lemma with \(I:=\Delta^2/(8\sigma^2)\). Hence
\[
\inf_{L_0\le L\le C/(1-p)}P_L\ge\frac{b\,(1-p)\log(1/(1-p))}{I}\,(1-o(1))=8(1-\delta)^2(u_g-\delta)\,v_g\,\kappa\log(1/\kappa)\,(1-o(1)),
\]
since \(\log(1/(1-p))=\log(1/\kappa)(1+o(1))\) for fixed \(\Delta/\sigma\). With the constant bound for \(L<L_0\), the same holds for the infimum over \(1\le L\le C/(1-p)\), and letting \(\delta\to0\) gives the lower bound \(8c_g\kappa\log(1/\kappa)(1-o(1))\), and the window of the upper bound has \(L=O(\log(1/\kappa))\le C/(1-p)\). The values of \(c_g\) follow from \(u_g=1/2\), \(v_g=1\) for uniform windows and \(u_g=\log2\), \(v_g=1/2\) for exponential weights.
\end{proof}

\begin{proof}[Proof of Proposition~\ref{prop:two-block}]
Let \(n:=a+b\) and \(\alpha:=\Delta^2/(8\sigma^2\|w\|^2)\). The \(\tau_t\) most recent observations come from the current regime, so \(W_t\) is at least the weight of the \(\tau_t\) most recent observations, which is \(k/(2a)\) for \(\tau_t=k\le a\), \(1/2+j/(2b)\) for \(\tau_t=a+j\) with \(0\le j\le b\), and \(1\) for \(\tau_t\ge n\). By Theorem~\ref{thm:margin} and \(\PP(\tau_t=k)\le1-p=\bar s\),
\[
\PP(\widehat E_t\neq E_t)\le\bar s\,a+\bar s\sum_{j=1}^{b}e^{-\alpha j^2/b^2}+e^{-\alpha}\le\bar s\,a+\bar s\,\frac b2\sqrt{\frac\pi\alpha}+e^{-\alpha},
\]
which holds for every \(\kappa\in(0,1)\). Since \(b\ge ra\), \(\|w\|^2=1/(4a)+1/(4b)\le(1+1/r)/(4a)\), so \(\alpha\ge\Delta^2a/(2\sigma^2(1+1/r))\ge\ell\) and \(e^{-\alpha}\le\kappa\). Moreover \(\bar s\,a\le2\kappa\ell(1+1/r)+\bar s\) and \(\bar s\,b\le r\bar s\,a+\bar s\), so
\[
\PP(\widehat E_t\neq E_t)\le2\kappa\ell\Bigl[1+\frac1r+(r+1)\frac{\sqrt\pi}{2\sqrt\ell}\Bigr]+\kappa+\bar s\Bigl(1+(r+1)\frac{\sqrt\pi}{2\sqrt\ell}\Bigr).
\]
With \(r=(4\ell/\pi)^{1/4}\), \(1/r+r\sqrt\pi/(2\sqrt\ell)=2(\pi/(4\ell))^{1/4}\), and since \(\bar s=\kappa\Delta^2/\sigma^2\) with \(\Delta/\sigma\) fixed, the remaining terms are \(O(\kappa\ell\cdot\ell^{-1/2})\). For Gaussian residuals the filter is optimal, so its error \eqref{eq:kz} bounds that of any decoder from below.
\end{proof}

\begin{proof}[Proof of Proposition~\ref{prop:centred}]
For the symmetric chain the switch indicators are independent Bernoulli\((1-p)\) variables, so at least two switches occur in the window with probability at most \(\binom{m-1}{2}(1-p)^2\le(1-p)^2(m-1)^2/2\). If no switch occurs, Theorem~\ref{thm:margin} with \(W_t=1\) bounds the error by \(e^{-m\Delta^2/(8\sigma^2)}\). With exactly one switch, at distance \(d\in\{0,\ldots,h-1\}\) from the centre on either side, the current regime occupies \(h+1+d\) positions, so \(2W_t-1=(2d+1)/m\), and Theorem~\ref{thm:margin} bounds the conditional error by \(\exp(-(2d+1)^2\Delta^2/(8\sigma^2m))\). Each of the \(2h\) positions has probability at most \(1-p\), and as in the proof of Corollary~\ref{cor:half-window},
\[
\sum_{d=0}^{h-1}\exp\Bigl(-\frac{(2d+1)^2\Delta^2}{8\sigma^2m}\Bigr)\le1+\frac12\int_0^\infty e^{-x^2\Delta^2/(8\sigma^2m)}\,dx=1+\frac{\sigma}{\Delta}\sqrt{\frac{\pi m}2}.
\]
This gives the middle term, and inserting the smallest odd \(m\) at least \((8\sigma^2/\Delta^2)\log(1/\kappa)\) gives the asymptotic form.
\end{proof}

\subsection{Proofs for Section~\ref{sec:k-regimes}}\label{app:proofs-k}

\begin{proof}[Proof of Proposition~\ref{prop:switch-geometry}]
\emph{Limits.} With \(k\) recent observations from regime~\(i\), the window average has mean \(\nu(k/m)\), where \(\nu(u):=\mu_j+u(\mu_i-\mu_j)\), and variance at most \(v/m\), where \(v\) is the larger of the two variances. As \(u\) increases from \(0\) to \(1\), \(\nu(u)\) moves from \(\mu_j\) to \(\mu_i\), and it lies in \(R_i\) exactly when \(u>c_{ji}\), and in \(R_j\) exactly when \(u<h_{ji}/|\mu_i-\mu_j|\). Fix \(\eta>0\) and let \(\varrho:=\min\{\eta|\mu_i-\mu_j|,\Delta/2\}\). If \(|k/m-c_{ji}|\ge\eta\), the distance from \(\nu(k/m)\) to the boundary of \(R_i\) is at least \(\varrho\), because the far end point of \(R_i\) lies beyond \(\mu_i\) at distance at least \(\Delta/2\). Chebyshev's inequality then gives
\[
\bigl|e_{ji}(k)-\mathbf 1\{k/m<c_{ji}\}\bigr|\le\frac{v}{m\varrho^2}.
\]
At most \(2\eta m+1\) values of \(k\) satisfy \(|k/m-c_{ji}|<\eta\), and the indicators average to \(c_{ji}+O(1/m)\), so \(\limsup_m|c_{ji}(m)-c_{ji}|\le2\eta\). Since \(\eta\) is arbitrary, \(c_{ji}(m)\to c_{ji}\). The same argument, applied to the event that the average lies outside \(R_i\cup R_j\), gives the limit of \(\chi_{ij}(m)\). Finally, \(h_{ij}\) is half the distance from \(\mu_i\) to its neighbour on the side of \(\mu_j\), so \(h_{ij}\le|\mu_i-\mu_j|/2\), with equality if and only if that neighbour is \(\mu_j\).

\emph{Gaussian expansion.} For Gaussian residuals, \(e_{ji}(k)=g(k)\) for the smooth function
\[
g(x):=\sum_{\theta}\Phi\Bigl(\pm\frac{(\theta-\nu(x/m))\sqrt m}{\sigma}\Bigr),
\qquad0\le x\le m,
\]
where the sum runs over the at most two end points \(\theta\) of \(R_i\), with the sign that makes each term the probability of falling beyond \(\theta\). The \(n\)-th derivative of \(g\) is \((|\mu_i-\mu_j|/(\sigma\sqrt m))^n\) times a derivative of \(\varphi\), so
\[
\int_0^m|g^{(n)}(x)|\,dx=O\bigl(m^{-(n-1)/2}\bigr).
\]
At \(x=0\) and \(x=m\) the mean \(\nu\) is \(\mu_j\) or \(\mu_i\), at distance at least \(\Delta/2\) from both end points, so every derivative of \(g\) at these points, as well as \(1-g(0)\) and \(g(m)\), is exponentially small in \(m\). For \(n\ge1\), the Euler--Maclaurin formula of order \(n\) has a remainder bounded by a constant multiple of \(\int_0^m|g^{(n)}(x)|\,dx\), so
\[
\sum_{k=1}^{m-1}g(k)=\int_0^mg(x)\,dx-\frac{g(0)+g(m)}{2}+O\bigl(m^{-(n-1)/2}\bigr)=\int_0^mg(x)\,dx-\frac12+O\bigl(m^{-(n-1)/2}\bigr).
\]
Since \(n\) is arbitrary, the error is \(O(m^{-r})\) for every \(r>0\). Finally, the term of \(g\) for the end point facing \(\mu_j\) integrates to \(m\,c_{ji}\) up to an exponentially small term, by the identity
\[
\int_{-a}^{b}\Phi(z)\,dz=b+\bigl[\varphi(b)-b\Phibar(b)\bigr]-\bigl[\varphi(a)-a\Phibar(a)\bigr],\qquad a,b>0,
\]
and the bounds \(0\le\varphi(x)-x\Phibar(x)\le\varphi(x)/x^2\) for \(x>0\), while the term for the far end point is exponentially small.
\end{proof}

\begin{proof}[Proof of Corollary~\ref{cor:k-window}]
(i) By Theorem~\ref{thm:switch-cost}, the error of \(m_{A_K}\) is at most \(A_K(m_{A_K})+\binom{m_{A_K}-1}2(1-p_{\min})\bar s\), and \(A_K(m_{A_K})\le A_K(m^\star)\le P^\star/p_{\min}^{m^\star-1}\).

(ii) Since \(Q=\mathrm{Id}+\varepsilon\Gamma\), the stationary law \(\pi\) and the ratios \(\bar s_{ji}/\bar s=\pi_j\Gamma_{ji}/\sum_{k\neq l}\pi_k\Gamma_{kl}\) do not depend on \(\varepsilon\), while \(\bar s\) and \(1-p_{\min}=\varepsilon\max_i|\Gamma_{ii}|\) are proportional to \(\varepsilon\). Write \(\ell:=\log(1/\bar s)\), \(a:=(1-p_{\min})/\bar s\), a constant, and \(P(m)\) for the error of the window of length \(m\). \emph{Long windows.} Let \(n_0:=\lceil1/\varepsilon\rceil\). For \(n\ge n_0\), \((Q^n)_{ji}=\sum_k(Q^{n-n_0})_{jk}(Q^{n_0})_{ki}\ge\min_k(Q^{n_0})_{ki}\), and \(Q^{n_0}\to e^{\Gamma}\), whose entries are positive because \(\Gamma\) is irreducible. The constant \(\rho\) of Lemma~\ref{lem:long-windows} depends only on \(\pi\) and \(\Gamma\), so for \(m-1\ge n_0\) and small \(\varepsilon\)
\[
P(m)\ge\frac{\rho}{1+\rho}\sum_{i\neq j,\ \Gamma_{ij}>0}\pi_j\,(Q^{m-1})_{ji}\ge c
\]
for a constant \(c>0\). By the upper bound below, \(P^\star\to0\), so for small \(\varepsilon\) the optimum lies among the windows of length at most \(2/\varepsilon\), which have length at most \(C'/\bar s\) for a constant \(C'\). We apply Lemma~\ref{lem:three-ranges} with \(s:=\bar s\) to these windows.

\emph{Single-regime errors.} As in the proof of Corollary~\ref{cor:window-rule}, Chernoff's bound and Cram\'er's theorem give, for every \(\delta>0\),
\[
\pi_{\min}\,c_\delta\,e^{-m(I+\delta)}\le\bar P(m)\le2e^{-mI},
\]
with a factor two because each decision region has at most two end points.

\emph{Costs.} Let
\[
r(m):=\frac1{\bar c\,\bar s}\sum_{i\neq j}\bar s_{ji}\,c_{ji}(m),
\qquad\text{so that}\qquad
A_K(m)=\bar P(m)+\bar c\,\bar s\,(m-1)\,r(m).
\]
Each \(c_{ji}(m)\) depends only on \(m\) and the emission laws, so \(r\) does not depend on \(\varepsilon\), and \(r(m)\to1\) by Proposition~\ref{prop:switch-geometry}. This is the uniformity over window lengths that the lemma requires.

\emph{Upper bound.} For \(m_0:=\lceil\ell/I\rceil\), Theorem~\ref{thm:switch-cost}(i) and~(ii) give
\[
P^\star\le A_K(m_0)+\frac{m_0^2}2\,(1-p_{\min})\,\bar s\le2\bar s+\frac{\bar c\,\bar s\,\ell}{I}\,r(m_0)+\frac a2\,\bar s^2m_0^2=\frac{\bar c\,\bar s\,\ell}{I}\,(1+o(1)).
\]

\emph{Lower bound.} By Theorem~\ref{thm:switch-cost}(ii), for every \(m\),
\[
P(m)\ge p_{\min}^{m-1}A_K(m)\ge(1-a\bar s)^m\bigl[\bar P(m)+\bar c\,\bar s\,(m-1)\,r(m)\bigr],
\]
which is condition~(iii) of the lemma with \(F:=\bar P\), \(b:=\bar c\) and \(\epsilon_s:=0\). The lemma gives \(P^\star\ge(\bar c\bar s\ell/I)(1-o(1))\) and \(m^\star=(\ell/I)(1+o(1))\). The function \(A_K\) satisfies the same condition and \(A_K(m_{A_K})\le A_K(m_0)\). Since \(r\ge1/2\) for large \(m\), \(A_K(m)\ge\bar c\,\bar s(m-1)/2\) there, so \(m_{A_K}=O(\ell)\), and the lemma gives \(m_{A_K}=(\ell/I)(1+o(1))\).

In the Gaussian case an end point at distance \(h\) from \(\mu_i\) has exponent \(h^2/(2\sigma^2)\), and the smallest \(h\) is \(\Delta/2\).
\end{proof}

\subsection{Proof for Section~\ref{sec:reflected}}\label{app:proofs-reflected}

\begin{proof}[Proof of Proposition~\ref{prop:reflected}]
Given \(\mathcal G_s\), \(Y_s\) takes values in an interval of length at most two, so Hoeffding's lemma gives \(\EE[e^{\lambda\eta_s}\mid\mathcal G_s]\le e^{\lambda^2/2}\). The argument of the proof of Proposition~\ref{thm:recoverability}, with \(\mathcal G_s\) in place of \(\mathcal H_s\) and \(\sigma=1\), gives
\[
\PP\Bigl(\Bigl|\sum_{s=t-m+1}^t\eta_s\Bigr|\ge\frac{m\Delta}{2}\Bigr)\le2e^{-m\Delta^2/8}.
\]
On \(\mathrm{NoSwitch}_{t,m}\cap H_{t,m}^c\), every step of the window starts in the interior under regime \(E_t\), so \(\bar Y_{t,m}-\mu_{E_t}=\frac1m\sum_s\eta_s\), and the decoder errs only if \(|\sum_s\eta_s|\ge m\Delta/2\).
\end{proof}

\clearpage
\section{Additional Experiments}\label{app:experiments}

\subsection{The half-window effect}

\begin{figure}[htbp]
\centering
\includegraphics[width=0.62\textwidth]{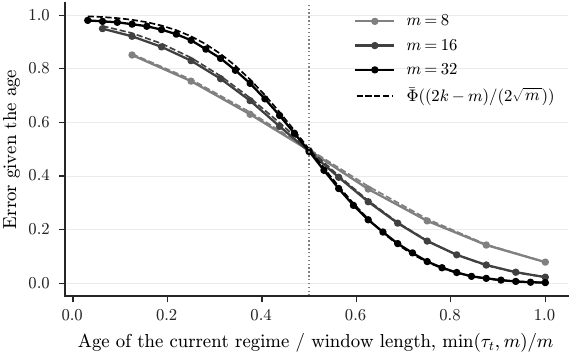}
\caption{Experiment~I, \(\Delta=\sigma=1\), \(p=\G{gauss-age-p}\). Error of the uniform window given the age of the current regime. Dashed: \(\Phibar\bigl((2k-m)/(2\sqrt m)\bigr)\).}
\label{fig:switch-age}
\end{figure}

Figure~\ref{fig:switch-age} isolates the half-window effect. For \(m=\G{gauss-age-m}\), a window in which the current regime is younger than \(m/2\) is misclassified with probability \(\G{gauss-age-before}\) on average, and one in which it is older with probability \(\G{gauss-age-after}\).

\subsection{\texorpdfstring{\(\pm1\)}{+-1} observations}

We compute the exact error over a sweep of \(\G{pm1-nconfig}\) pairs \((p,\Delta)\) and check it by simulation (\(\G{pm1-R}\) trajectories of length \(\G{pm1-T}\)). The two agree to within \(\G{pm1-maxdev}\), or \(\G{pm1-maxz}\) standard errors, and the simulated optimum coincides with the exact one in \(\G{pm1-nagree-sim}\) configurations. Here windows are short, with \(m^\star\) ranging from \(\G{pm1-mstar-lo}\) to \(\G{pm1-mstar-hi}\), and the window of length one is optimal in \(\G{pm1-nones}\) configurations. Every optimal window is odd. Ties explain this.

\begin{lemma}[Even windows]
\label{lem:even-windows}
Under model~(a) with ties broken by a fair coin, \(\bar P(2k)=\bar P(2k-1)\) for every \(k\ge1\).
\end{lemma}

\begin{proof}
In a single-regime window each observation points to the wrong regime independently with probability \(q:=(1-\Delta/2)/2\). Let \(W\) be the number of wrong observations among the first \(2k-1\). Then \(\bar P(2k-1)=\PP(W\ge k)\), and
\[
\bar P(2k)=\PP(W\ge k+1)+\PP(W=k)\Bigl(q+\frac{1-q}{2}\Bigr)+\PP(W=k-1)\,\frac q2 ,
\]
which differs from \(\PP(W\ge k)\) by \(\bigl[q\,\PP(W=k-1)-(1-q)\,\PP(W=k)\bigr]/2=0\).
\end{proof}

\begin{figure}[htbp]
\centering
\includegraphics[width=\textwidth]{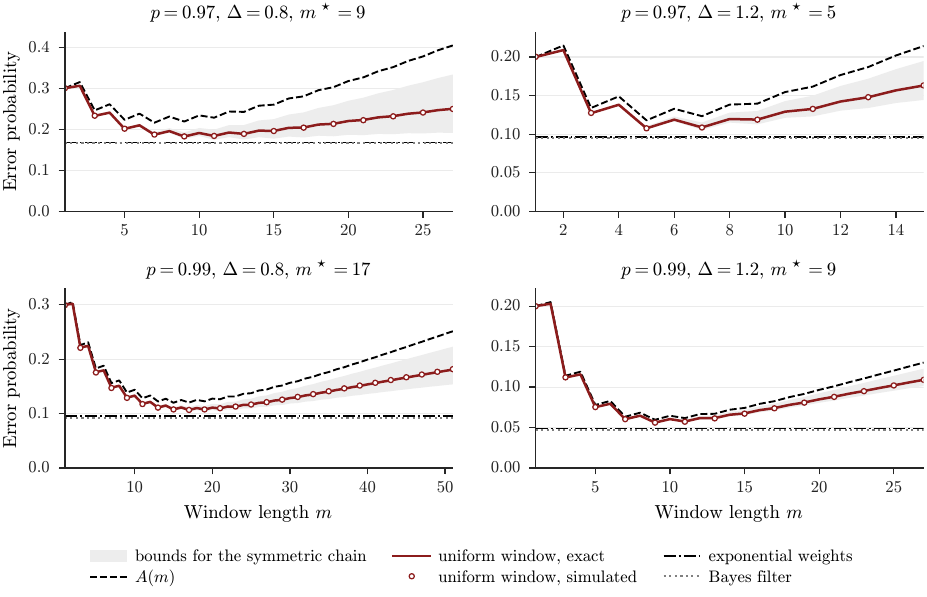}
\caption{Experiment~I, \(\pm1\) observations. Error of the uniform window for windows up to three times the optimum.}
\label{fig:expII-window-tradeoff}
\end{figure}

An even window therefore has the single-regime error of the preceding odd one and a larger chance of containing a switch. At \(\Delta=\G{pm1-tie-delta}\) a window of length two ties with probability \(\G{pm1-tie-prob}\). The minimizer \(m_A\) of the half-window law is the optimal window in \(\G{pm1-nA-exact}\) of the \(\G{pm1-nconfig}\) configurations, with a median excess error of \(\G{pm1-A-excess-med}\%\) and a maximum of \(\G{pm1-A-excess-max}\%\), while the minimizer of \(g_2\) costs up to \(\G{pm1-g2-excess-hi}\%\). Table~\ref{tab:expII} and Figure~\ref{fig:expII-window-tradeoff} show four configurations in detail. In them the bounds of Corollary~\ref{cor:symmetric-bounds} enclose the exact error in a band that is narrow near the optimum, exponential weights are within \(\G{pm1-e-filt-lo}\)--\(\G{pm1-e-filt-hi}\) times the filter error, and the uniform window within \(\G{pm1-u-filt-lo}\)--\(\G{pm1-u-filt-hi}\).

\begin{table}[htbp]
\centering
\caption{Experiment~I, \(\pm1\) observations. Optimal window \(m^\star\), the minimizer \(m_A\) of the half-window law and their errors (exact), with exponential weights and the Bayes filter (simulated, standard errors in parentheses).}
\label{tab:expII}
\begin{tabular}{cccccccc}
\toprule
\(p\) & \(\Delta\) & \(m^\star\) & \(m_A\) & error at \(m^\star\) & error at \(m_A\) & exponential & Bayes filter \\
\midrule
0.97 & 0.80 & 9 & 7 & 0.1823 & 0.1868 & 0.1676\,(0.0004) & 0.1660\,(0.0004) \\
0.97 & 1.20 & 5 & 5 & 0.1076 & 0.1076 & 0.0967\,(0.0003) & 0.0955\,(0.0003) \\
0.99 & 0.80 & 17 & 15 & 0.1066 & 0.1076 & 0.0948\,(0.0004) & 0.0918\,(0.0005) \\
0.99 & 1.20 & 9 & 9 & 0.0561 & 0.0561 & 0.0488\,(0.0002) & 0.0469\,(0.0002) \\
\bottomrule
\end{tabular}

\end{table}

\FloatBarrier

\authorcontact
  {Arina Pambukyan}
  {American University of Armenia, Yerevan, Armenia}
  {arina\string_pambukyan26@alumni.aua.am}

\authorcontact
  {Hayk Nersisyan}
  {American University of Armenia, Yerevan, Armenia}
  {hnersisyan@aua.am}


\begin{thebibliography}{99}

\bibitem{AbrahamGassiatNaulet2023}
K. Abraham, E. Gassiat and Z. Naulet.
\newblock Fundamental limits for learning hidden Markov model parameters.
\newblock {\em IEEE Transactions on Information Theory}, 69(3):1777\textendash 1794, 2023.

\bibitem{AngTimmermann2012}
A. Ang and A. Timmermann.
\newblock Regime changes and financial markets.
\newblock {\em Annual Review of Financial Economics}, 4:313\textendash 337, 2012.

\bibitem{Asmussen2003}
S. Asmussen.
\newblock {\em Applied Probability and Queues}.
\newblock 2nd edition, Springer, 2003.

\bibitem{Azuma1967}
K. Azuma.
\newblock Weighted sums of certain dependent random variables.
\newblock {\em T\^ohoku Mathematical Journal}, 19(3):357\textendash 367, 1967.

\bibitem{BassevilleNikiforov1993}
M. Basseville and I. V. Nikiforov.
\newblock {\em Detection of Abrupt Changes: Theory and Application}.
\newblock Prentice Hall, 1993.

\bibitem{BootPick2018}
T. Boot and A. Pick.
\newblock Optimal forecasts from Markov switching models.
\newblock {\em Journal of Business \& Economic Statistics}, 36(4):628\textendash 642, 2018.

\bibitem{CappeMoulinesRyden2005}
O. Capp\'e, E. Moulines and T. Ryd\'en.
\newblock {\em Inference in Hidden Markov Models}.
\newblock Springer, 2005.

\bibitem{Chigansky2006}
P. Chigansky.
\newblock On filtering of Markov chains in strong noise.
\newblock {\em IEEE Transactions on Information Theory}, 52(9):4267\textendash 4272, 2006.

\bibitem{ChoiYeungZhang2000}
S. P. M. Choi, D.-Y. Yeung and N. L. Zhang.
\newblock An environment model for nonstationary reinforcement learning.
\newblock In {\em Advances in Neural Information Processing Systems 12}, pages 987\textendash 993, 2000.

\bibitem{DaiZhangZhu2010}
M. Dai, Q. Zhang and Q. J. Zhu.
\newblock Trend following trading under a regime switching model.
\newblock {\em SIAM Journal on Financial Mathematics}, 1(1):780\textendash 810, 2010.

\bibitem{DemboZeitouni1998}
A. Dembo and O. Zeitouni.
\newblock {\em Large Deviations Techniques and Applications}.
\newblock 2nd edition, Springer, 1998.

\bibitem{EichingerKirch2018}
B. Eichinger and C. Kirch.
\newblock A MOSUM procedure for the estimation of multiple random change points.
\newblock {\em Bernoulli}, 24(1):526\textendash 564, 2018.

\bibitem{GajaneOrtnerAuer2018}
P. Gajane, R. Ortner and P. Auer.
\newblock A sliding-window algorithm for Markov decision processes with arbitrarily changing rewards and transitions.
\newblock arXiv:1805.10066, 2018.

\bibitem{GarivierMoulines2011}
A. Garivier and E. Moulines.
\newblock On upper-confidence bound policies for switching bandit problems.
\newblock In {\em Algorithmic Learning Theory}, Lecture Notes in Computer Science 6925, pages 174\textendash 188. Springer, 2011.

\bibitem{GassiatKaddouriNaulet2025}
E. Gassiat, I. Kaddouri and Z. Naulet.
\newblock Clustering risk in nonparametric hidden Markov and i.i.d.\ models.
\newblock {\em The Annals of Statistics}, 53(6):2409\textendash 2429, 2025.
\newblock \href{https://doi.org/10.1214/25-AOS2546}{doi:10.1214/\allowbreak\mbox{25-AOS2546}}.

\bibitem{Hamilton1989}
J. D. Hamilton.
\newblock A new approach to the economic analysis of nonstationary time series and the business cycle.
\newblock {\em Econometrica}, 57(2):357\textendash 384, 1989.

\bibitem{KhammassiEtAl2025}
M. Khammassi, V. Bordignon, V. Matta and A. H. Sayed.
\newblock Fundamental social learning scaling law for tracking hidden {M}arkov models.
\newblock In {\em Proceedings of the IEEE International Conference on Acoustics, Speech and Signal Processing (ICASSP)}, 2025. \href{https://doi.org/10.1109/ICASSP49660.2025.10888523}{doi:10.1109/\allowbreak\mbox{ICASSP49660.2025.10888523}}.

\bibitem{KhasminskiiZeitouni1996}
R. Khasminskii and O. Zeitouni.
\newblock Asymptotic filtering for finite state Markov chains.
\newblock {\em Stochastic Processes and their Applications}, 63(1):1\textendash 10, 1996.

\bibitem{KirchKlein2023}
C. Kirch and P. Klein.
\newblock Moving sum data segmentation for stochastic processes based on invariance.
\newblock {\em Statistica Sinica}, 33(2):873\textendash 892, 2023.

\bibitem{Lai1974}
T. L. Lai.
\newblock Control charts based on weighted sums.
\newblock {\em The Annals of Statistics}, 2(1):134\textendash 147, 1974.

\bibitem{LangrockEtAl2012}
R. Langrock, R. King, J. Matthiopoulos, L. Thomas, D. Fortin and J. M. Morales.
\newblock Flexible and practical modeling of animal telemetry data: hidden Markov models and extensions.
\newblock {\em Ecology}, 93(11):2336\textendash 2342, 2012.

\bibitem{LundeTimmermann2004}
A. Lunde and A. Timmermann.
\newblock Duration dependence in stock prices: an analysis of bull and bear markets.
\newblock {\em Journal of Business \& Economic Statistics}, 22(3):253\textendash 273, 2004.

\bibitem{PambukyanNersisyanCode2026}
A. Pambukyan and H. Nersisyan.
\newblock Sliding-window detection of hidden regimes: code, version 1.0.0.
\newblock Zenodo, 2026.
\newblock \href{https://doi.org/10.5281/zenodo.24532266}{doi:10.5281/\allowbreak\mbox{zenodo.24532266}}.

\bibitem{PambukyanSaudagarKumar2026}
A. Pambukyan, A. K. J. Saudagar and S. Kumar.
\newblock Ergodicity and regime recoverability in finite Markov-modulated random walks.
\newblock {\em Scientific Reports}, 16:12376, 2026.
\newblock \href{https://doi.org/10.1038/s41598-026-45019-w}{doi:10.1038/\allowbreak\mbox{s41598-026-45019-w}}.

\bibitem{PerssonEtAl2013}
F. Persson, M. Lind\'en, C. Unoson and J. Elf.
\newblock Extracting intracellular diffusive states and transition rates from single-molecule tracking data.
\newblock {\em Nature Methods}, 10(3):265\textendash 269, 2013.

\bibitem{PrabhuZhu1989}
N. U. Prabhu and Y. Zhu.
\newblock Markov-modulated queueing systems.
\newblock {\em Queueing Systems}, 5:215\textendash 246, 1989.

\bibitem{SchwierBrooksGriffin2011}
J. M. Schwier, R. R. Brooks and C. Griffin.
\newblock Methods to window data to differentiate between Markov models.
\newblock {\em IEEE Transactions on Systems, Man, and Cybernetics, Part B (Cybernetics)}, 41(3):650\textendash 663, 2011.

\bibitem{SheikhBhat2026}
S. A. Sheikh and D. A. Bhat.
\newblock Pockets of predictability: a stochastic crowding model and finite-sample regime detection for the endogenous decay of trading anomalies.
\newblock {\em Computational Mathematics and Modeling}, 2026. \href{https://doi.org/10.1007/s10598-026-09751-2}{doi:10.1007/\allowbreak\mbox{s10598-026-09751-2}}.

\bibitem{XieMoustakidesXie2023}
L. Xie, G. V. Moustakides and Y. Xie.
\newblock Window-limited {CUSUM} for sequential change detection.
\newblock {\em IEEE Transactions on Information Theory}, 69(9):5990\textendash 6005, 2023.

\bibitem{ZhaoEtAl2026}
Y. Zhao, O. Eberhard, M. Khammassi, A. H. Sayed and M. Muehlebach.
\newblock Why linear recurrent memory works in partially observable reinforcement learning.
\newblock In {\em Proceedings of the 43rd International Conference on Machine Learning}, PMLR 306, pages 161337\textendash 161372, 2026.

\bibitem{ZhouXiongChenGao2021}
X. Zhou, Y. Xiong, N. Chen and X. Gao.
\newblock Regime switching bandits.
\newblock In {\em Advances in Neural Information Processing Systems 34}, 2021.

\end{thebibliography}
\end{document}